\documentclass{article}

\usepackage{amsmath,amssymb,amsthm,float,mathrsfs}
\usepackage[dvipdfmx]{graphicx,xcolor}
\usepackage{comment}
\usepackage{authblk}

\newcommand{\1}[1]{1_{\{#1\}}}
\newcommand{\norm}[1]{\lVert  #1 \rVert}
\newcommand{\C}{\mathcal{C}}

\newcommand{\E}{\mathbb{E}}
\newcommand{\HP}{\widehat{\mathbb{P}}}
\newcommand{\HE}{\widehat{\mathbb{E}}}
\newcommand{\CE}{\mathcal{E}}
\newcommand{\CF}{\mathcal{F}}
\newcommand{\CH}{\mathcal{H}}
\renewcommand{\P}{\mathbb{P}}

\newcommand{\hl}{\widehat{\mathbb{L}}^2}

\newcommand{\setmid}{\mathrel{}\middle|\mathrel{}}
\renewcommand{\epsilon}{\varepsilon}
\newcommand{\OH}{\widehat{\Omega}}
\newcommand{\cp}[1]{\text{Cap}\left(#1\right)}
\renewcommand{\phi}{\varphi}

\theoremstyle{definition}
\newtheorem{Def}{Definition}[section]
\newtheorem{thm}[Def]{Theorem}
\newtheorem{lem}[Def]{Lemma}

\newtheorem{prop}[Def]{Proposition}
\newtheorem{cor}[Def]{Corollary}
\newtheorem{rem}[Def]{Remark}

\newtheorem{eg}[Def]{Example}

\newtheorem{asm}{Assumption}

\makeatletter

\@addtoreset{equation}{section}
\makeatother

\begin{document}
	\title{Quenched Invariance Principle for a Reflecting Diffusion in a Continuum Percolation Cluster}
	\author{Yutaka \textsc{takeuchi} \thanks{Department of Mathematics, Faculty of Science and Technology, Keio University \\ e-mail: \texttt{yutaka.takeuchi@keio.jp}}}
	\date{}
    \maketitle
	\section{Introduction}
	\subsection{Models and result}
	Let $ \Omega $ be a configuration space on $ \mathbb{R}^d $, namely \[\Omega = \left \{\sum _{i=1}^\infty \delta _{x_i}  \setmid \ \sum _{i=1}^\infty \delta _{x_i}(K) <\infty \text{ for all compact set }K  \right \} .\]
	Here $\delta_x$ denotes the Dirac measure.
	$ \Omega $ is equipped with the $ \sigma $-field $ \mathcal{B}(\Omega) $ which is generated by the sets $ \{ \omega \in \Omega \mid \omega(A) = n\} $, $ A \in \mathcal{B}(\mathbb{R}^d) , n \in \mathbb{N}$.  
	We identify each configuration $\omega = \sum _{i=1}^\infty \delta _{x_i} \in \Omega $ with the subset $ \{x_i\}_i $ of $ \mathbb{R}^d $ if $ \{x_i\}_i $ are distinct. 
	Fix $\rho>0$. For $ \omega \in \Omega $, define the subset $ L(\omega)$ by
	\[ L(\omega) =\bigcup _{x\in \omega}B(x,\rho),\]
	where $B(x,\rho)$ is the Euclidean open ball with center $x$ and radius $\rho$. Let $W(\omega)$ be the unbounded connected component of $L(\omega)$ if there is a unique unbounded component. Otherwise, we set $ W(\omega)=\emptyset $ by convention. We call $W(\omega)$ the continuum percolation cluster. This is the continuum analogue of the discrete site percolation cluster. 
	By definition, $W(\omega)$ can be written as $\bigcup_{x \in I(\omega)}B(x,\rho)$ with some subset $I(\omega)$ of $\omega$. We also take $\rho' \geq \rho$ and we introduce the modified cluster $W'(\omega)$ given by
	\begin{align*}
		W'(\omega) \equiv W_{\rho'}(\omega) = \bigcup_{x \in I(\omega)}B(x,\rho'). 
	\end{align*}
	Set 
	\[\Delta = \left \{\omega  \in \Omega \setmid \begin{aligned}
	&\omega = \sum _i\delta _{x_i},\;\;  |x-y|=2\rho \;\;  (\text{for some } x,y \in \omega ) \\ &\text{ and each $x_i$ is distinct}
	\end{aligned}   \right \}.\] 

	Define the subset $ \OH $ of $ \Omega $ by
	\[\OH=\{\omega \in \Omega \setminus \Delta \mid  0 \in  W'(\omega) \}.\]

	We define the shift  $\tau _z\colon \Omega \to \Omega$ by
	\[\tau _z \omega (A) = \omega (A+z)=\sum_{x\in \omega} \delta_{x}(A+z)=\sum_{x\in \omega} \delta_{x-z}(A),\]
	where $A+z = \{x+z \mid x \in A\}$ for $z \in \mathbb{R}^d$. 
	
	Let $\P $ be a probability measure on $\Omega$. It is called a point process on $ \mathbb{R}^d $. 
	\begin{asm}\label{asm:erg} Assume that $ \P $ satisfies the following conditions:
		\begin{enumerate} 
			\item $\P$ is stationary and ergodic with respect to $\{\tau _x\}_x$.
			\item $\P (\OH)>0$ and $\P(\Delta) = 0$. 
		\end{enumerate}
	\end{asm}
	Set $\HP (\cdot )=\P (\cdot  \cap  \OH)/\P(\OH)$ and we denote its expectation by $ \HE $.
	Note that $W'(\omega)$ is a Lipschitz domain for $\HP$-a.e. $\omega \in \OH$. Thanks to $\P(\Delta) = 0$, we only need to consider whether $ 0 \in W(\omega) $ or not. The existence of a unique unbounded component is an important problem in the study of continuum percolation theory and there are many studies. If $ \P $ is a Poisson point process and the radius $ \rho $ is bigger than the critical value, it is known that there is a unique unbounded component (see \cite{MR}). 
	
	We denote the Euclidean inner product by $ \langle \cdot , \cdot \rangle $. 
	Let $ a\colon \Omega \to \mathbb{R}^{d\times d} $ be a positive-definite symmetric matrix. Define the bilinear form $ \CE $ on $L^2(W'(\omega),dx)$ by
	\begin{align}
	\CE^\omega (u,v)=\int_{W'(\omega)} \langle a(\tau_x\omega) \nabla u(x), \nabla v(x) \rangle dx.
	\end{align}
	Let $\CF^\omega$ be the completion of $C_c^\infty(\overline{W})$ with respect to $\CE(\cdot,\cdot)+\norm{\cdot}_{L^2(W'(\omega))}$.
	\begin{asm}\label{asm:ellipse}
		There exist constants $\lambda,\Lambda > 0$ such that
		\begin{align*}
			\lambda |\xi|^2 \leq \langle a(\omega)\xi, \xi \rangle \leq \Lambda |\xi|^2
		\end{align*}
		holds for all $\xi \in \mathbb{R}^d$ and $\HP$-almost all $\omega$.
	\end{asm}

	According to \cite{FT} and \cite{FT2}, under Assumption \ref{asm:ellipse}, we have that the Dirichlet form $(\CE^\omega,\, \CF^\omega)$ is strongly local and regular $\HP$-almost surely.  Hence, we have the associated conservative diffusion $X_t^\omega$. It is called a reflecting diffusion since the domain $\CF^\omega$ of Dirichlet form corresponds to a reflecting boundary condition (more precisely, the Neumann boundary condition). 
	
	We further impose an assumption for reflecting diffusions:
	\begin{asm}\label{asm:density}
		The reflecting diffusion $X_t^\omega$ has a transition density $p_t^\omega(\cdot,\cdot)$ for $\HP$-almost all $\omega$.
	\end{asm}
	When we consider the case that $a(\omega)=1/2I_d$, $X_t$ is the reflecting Brownian motion. Bass and Hsu (\cite{BH}) study a condition under which $X_t$ has a density. Recently, Matsuura (\cite{Ma}) considered more detailed conditions. Our model satisfies the conditions of Matsuura \cite{Ma}, and hence the reflecting Brownian motion has a density.

	$ $

	Next, we impose a geometric condition that plays an important role. 
	Let $ x \in \mathbb{R}^d $ and $ R > 0 $. 	
	Let $W'(\omega, x, R)$ be a connected component of $W'(\omega)\cap B(x,R)$ containing $x$. Throughout the paper, we write $ W'=W'(\omega) $ and $ W'_R $ = $ W'(\omega, 0, R) $ if there is no confusion. 
	
	\begin{asm}[volume regularity and isoperimetric condition]\label{asm:volIso}

		$ $

		\begin{enumerate}
			\item For $\HP$-almost every $ \omega $ there exists positive constant $C_V$ such that for a.e. $x\in W'(\omega)$ 
        \[ C_V R^d \leq |W'(\omega, x, R)|,  \quad R\geq R_V\]
        holds for some positive constant $R_V$, where $ |\cdot| $ denotes the Lebesgue measure. 	
		\item There exists $\theta \in (0,1)$ such that for $\HP$-almost all $\omega \in \OH$, there exist $R_I>0$ and $c_H>0$ such that
        \begin{gather}\label{asm:IPforLarge}
			C_\text{IL} := \inf \left\{\frac{\CH_{d-1}(W' \cap \partial O)}{|W' \cap O |^\frac{d-1}{d}} \middle | \begin{split}
				&\text{ $O\subset B(0,R)$ is connected open,}\\&R\geq R_I,  \;|O| \geq R^\theta  \end{split}  \right\} > 0
		\end{gather}
		and
		\begin{gather}\label{asm:IPforShort}
			C_\text{IS}:=\inf\left\{\frac{\mathcal{H}_{d-1}(W'\cap \partial O)}{|W'\cap O|^\frac{d-1}{d}}\middle | \begin{split} &O\text{ is bounded open,}\\ &\mathcal{H}_{d-1}(W'\cap \partial O) < c_H
			\end{split}\right\} > 0
		\end{gather}
		hold, where $\mathcal{H}_{d-1}$ denotes the $(d-1)$-dimensional Hausdorff measure.
		\end{enumerate}
	\end{asm}
	
	The main result of this paper is the following.
	
	\begin{thm}[quenched invariance principle]
		Let $ d \geq 2 $ and $\rho' \geq \rho> 0$. Assume that Assumptions \ref{asm:erg}, \ref{asm:ellipse}, \ref{asm:density} and \ref{asm:volIso} hold. Let $ P_0^\omega $ be the law of $ \{X_t^\omega\}_t$ starting at $ 0 $. \\
		Then for $\HP$-a.s.$\omega$, the scaled process $\{\epsilon X_{\epsilon ^{-2}t}^\omega\}_t$ under $ P_0^\omega $ converges in law to a Brownian motion with non-degenerate covariance matrix as $\epsilon$ tends to $0$.
	\end{thm}
	
\begin{rem}\label{rem:thm}
	 When $\rho' > \rho$, the condition $(\ref{asm:IPforShort})$ is automatically satisfied. On the other hand, the condition $(\ref{asm:IPforShort})$ doesn't hold when $\rho'=\rho$ in general. 
\end{rem}

	\begin{eg}[Reflecting Brownian motion on Poisson Boolean model]\label{eg:PBM}
		Let $\P$ be a Poisson point process with intensity $\lambda$. We take the radius $\rho$ greater than the critical radius $\rho_c$. We also take the radius $\rho'>\rho$. Set $a(\omega) = 1/2I_d$. Then  the corresponding diffusion is the reflecting Brownian motion. We can easily verify that Assumption \ref{asm:erg} and \ref{asm:ellipse} are satisfied. Assumption \ref{asm:density} is satisfied by the result in \cite{Ma}. Now we need to check the Assumption \ref{asm:volIso}. We check this by comparing the modified cluster $W'(\omega)$ to the Bernoulli site percolation model as in \cite{T2}.
		To do this, we prepare a percolation model on $\delta \mathbb{Z}^d$, $\delta >0$. We call the model the $\delta$-approximating Boolean model. 
		Set $G(z,\delta) = [-\delta/2,\delta/2]^d+z$. Take $\delta \in (0,1) $. We say that a site $z \in \delta\mathbb{Z}^d$ is $\delta$-open if $G(z,\delta) \subset W'(\omega)$.
		Let $V_\delta$ be the collection of $\delta$-open sites on $\delta \mathbb{Z}^d$. 
		We introduce the graph $\mathcal{G}_\delta = (V_\delta, E_\delta)$ with vertex set $V_\delta$ and edge set $E_\delta =\{\{z,z'\} \subset V_\delta\mid |z-z'|=\delta\}$. For $x \in V_\delta$, let $B_\delta^\omega(x,R)$ be an open ball of $V_\delta$ with respect to the graph distance.
		For a cube $Q = G(z,n\delta)$, write $Q^+=G(z,\frac{3}{2}n\delta)$. We say that a cluster $\mathcal{C}$ in a cube $Q$ is crossing for a cube $Q'\subset Q$ if for all $d$-directions there exists an open path in $\mathcal{C} \cap Q'$ that connects the two opposing faces of $Q'$. Let $\C^\vee(Q)$ be the largest cluster in $Q$. As in \cite{B}, set
		\begin{align*}
			\begin{split}
			R_0(Q)=\{ 
					&\text{there exists a unique crossing cluster } \mathcal{C} \text{ in }Q^+ \text{ for } Q^+  \text{ such that} \\
					&\text{ all open paths contained in } Q^+ \text{ of diameter greater than } n\delta/8\\
					& \text{ are connected to } \mathcal{C} \text{ in } Q^+ \text{ and } \mathcal{C} \text{ is crossing for each cube } \\
					&\text{contained in } Q \text{ whose side length is greater than or equal to } n\delta/8\}
				\end{split}
		\end{align*}
		and 
		\begin{align*}
			R(Q) = R_0(Q)\cap\{\C^\vee(Q)\text{ is crossing for }Q\} \cap \{\C^\vee( Q^+)\text{ is crossing for }Q^+\}.
		\end{align*}

		First we consider the following estimate:
		\begin{align}\label{ine:Rexp}
			\P(R(G(z,k\delta))^c)\leq c \exp(-c'k), \; k\in \mathbb{N}.
		\end{align}
		In \cite[Lemma 2.8]{B}, the estimate $(\ref{ine:Rexp})$ is proved in the case of the Bernoulli site percolation. We can generalize the estimate $(\ref{ine:Rexp})$ to the $\delta$-approximating Boolean model by renormalization used in \cite{T2}. 
		Hence as in \cite{B}, we can show that the process $\{1_{R(G(z,k\delta))}\}_{z\in\delta\mathbb{Z}^d}$ dominates the Bernoulli site percolation with parameter $q_\delta^*(k)$, which tends to $1$ as $k \to \infty$ and satisfies $q_\delta^*(k) \leq q_{\delta'}^*(k)$ for $\delta \geq \delta'$.  
		Therefore, for sufficiently large $n$ and a cube $Q=G(z_0, n\delta)$, we can prove the weak relative isoperimetric inequality 
		\begin{align}
			\#\{ \{z, z'\}  \mid z\in A,\; z'\in \mathcal{C}^\vee(Q)-A,\; |z-z'|=\delta  \} \geq c_1 n^{-1}\#A
		\end{align}
		 holds for connected subset $A \subset \mathcal{C}^\vee(Q)$ with $\#A\leq 1/2\#\mathcal{C}^\vee(Q)$ such that $\mathcal{C}^\vee(Q)-A$ is also connected, and the volume regularity
		\begin{align}\label{ine:appoxVolReg}
			c_2R^d\leq \#B_\delta^\omega(z,R) \leq c_2'R^d
		\end{align}
		holds for all $z \in \mathcal{C}^\vee(Q^+)\cap G(z_0, 5/6n\delta)$ with $G(z,R+k)^+ \subset Q^+$ and $R \in (c_Hn^\alpha, n)$ as in \cite[Proposition 2.11 and Theorem 2.23]{B}.	
		Then by \cite[Lemma 2.10]{DNS}, there exist $\theta \in (0,1)$ and $N_3$ such that  $V_\delta$ satisfies the isoperimetric inequality for large sets
		\begin{align}\label{ine:approxIso}
			\frac{\# \{\{z,z'\} \in E_\delta \mid z\in A,  z' \in V_\delta - A\}}{(\# A)^\frac{d-1}{d}}\geq c_3 
		\end{align}
		holds for $R \geq N_3$ and $A \subset B_\delta^\omega(0, R)$ with $\#A \geq R^\theta$. Note that $B_\delta^\omega(0,R)$ is contained in the unique infinite cluster of $V_\delta$. 
		Furthermore, we can take constants $c_1, c_2, c_2', c_3 > 0$, and $N_3>0$ 
		independently of $\delta \in (0,\delta_0)$ for some $\delta_0$ 
		(This is because these constants depend only on $k$ and $q_\delta^*(k)$ is decreasing in $\delta$.).

		Now we can check Assumption \ref{asm:volIso}. Since $\rho'>\rho$, the condition $(1.3)$ is satisfied (see Remark \ref{rem:thm}.).  Because $B_\delta^\omega(x,R/\delta)$ is contained in $W'(\omega,x,R)$, we have $|W'(\omega,x,R)| \geq \delta^d \# |B_\delta^\omega(x,R/\delta)|$. Hence we can easily see that (1) of Assumption \ref{asm:volIso} holds.
		To check (2) of Assumption \ref{asm:volIso}, we introduce some notations. Set 
		\begin{align*}
			&\mathbb{G}_\delta = \{G(z,\delta) \mid z \in \mathbb{Q}^d\},\; \;\mathbb{G} = \bigcup_{\delta \in (0,1) \cap \mathbb{Q}} \mathbb{G}_\delta, \;\;
			\mathcal{Q}_\text{fin}=\left\{\bigcup_{i=1}^N Q_i \mid Q_i \in \mathbb{G} \right\}.  
		\end{align*}
		Take a connected subset $D  \in \mathcal{Q}_\text{fin}$ with $W'(\omega) \cap D \neq \emptyset$.
		  We can choose $\delta >0$ so that $D$ is written as  $\bigcup_{i=1}^N G(z_i,\delta)$, $z_i\in \delta\mathbb{Z}^d$, $i=1,\dots,N$. Moreover, we can take $\delta>0$ smaller if we need. Set $Z =\{z_1, \dots, z_N\}$. Let $Z_o$ be the collection of $\delta$-open sites of $Z$ and $\partial_\text{ext} Z_o$ be the external boundary of $Z_o$.
		Observe that at least one face of cubes in $\bigcup_{z\in \partial_\text{ext}Z_o}G(z,\delta)$ belongs to  $W'(\omega) \cap \partial D$. Therefore we have $\mathcal{H}_{d-1}(W'(\omega) \cap \partial Q) \geq \delta^{d-1} \# \partial_\text{ext}Z_o$. Next, let $Z_c = \{z \in Z\setminus Z_o \mid G(z,\delta)\cap W'(\omega) \neq \emptyset\}$. Then  we have that $|W'(\omega)\cap D| \leq \delta^d(\#Z_o + \#Z_c)\leq 2\delta^d\#Z_o $ for sufficiently small $\delta$.
		Hence if $|D|\geq R^\theta$, using the isoperimetric inequality $(\ref{ine:approxIso})$,  we have that $\mathcal{H}_{d-1}(W'(\omega)\cap \partial D)\geq C_1 |W'(\omega)\cap D|^\frac{d-1}{d}$.  
		For a general connected open subset $O \subset B(0,R)$, we can take a sequence of $D_i \in \mathcal{Q}_\text{fin}$ such that $|W'(\omega)\cap D_i| \to |W'(\omega) \cap O|$ and $\mathcal{H}_{d-1}(W'(\omega)\cap \partial D_i) \to \mathcal{H}_{d-1}(W'(\omega)\cap \partial O)$. Hence, the condition $(\ref{asm:IPforLarge})$ holds.

		\begin{eg}
			Let $\{Z_z\}_{z\in \mathbb{Z}^d}$ be a collection of i.i.d. Bernoulli random variables with parameter $p \in (0,1)$ and $U$ be the uniform random variable on $[0,1)^d$. Let $\P$ be the distribution of $\sum_{z\in\mathbb{Z}}Z_z\delta_{z+U}$. 
			Let $\rho_c > 0$ be the critical radius. We remark that $\rho_c = 1/2$ if the parameter $p$ is greater than the critical value of the Bernoulli site percolation. Take $\rho > \rho_c$ and let $\rho'=\rho$. Set $a(\omega) = 1/2I_d$. Then we can check Assumption \ref{asm:erg}, \ref{asm:ellipse}, and \ref{asm:density} hold. Similarly to Example 1.3, we can check (1) of Assumption \ref{asm:volIso} and $(\ref{asm:IPforLarge})$ hold.
			To check $(\ref{asm:IPforShort})$, observe that by construction, there is only a finite number of possible distances between overlapping balls. Indeed, their values are of the form $|x|$, $x \in \mathbb{Z}^d$, and the maximum is strictly less than $2\rho$.  Therefore, computing $\frac{\mathcal{H}_{d-1}(W'\cap \partial O)}{|W'\cap O|^\frac{d-1}{d}}$ in each case, we can check $(\ref{asm:IPforShort})$.
			\end{eg}

	\end{eg}

	Early in the history of homogenization, Kipnis and Varadhan\cite{KV} proved that the annealed invariance principle for random walk on supercritical (bond) percolation cluster. After two decades, Sidoravicius and Sznitman \cite{SS} proved the quenched invariance principle for random walk on supercritical percolation cluster when the dimension is more than or equal to four. In 2007, Berger and Biskup \cite{BB} proved quenched invariance principle for random walk on the supercritical percolation cluster including two and three dimensions. (Note that one-dimensional percolation cluster is infinite only when the probability that bond is open is equal to one and in that case the QIP follows from classical Donsker's invariance principle.) Mathieu and Piatnitski \cite{MP} gave another proof. Quenched invariance principle for the random conductance model, which is a more general model, was shown by many authors (\cite{A}, \cite{ABDH}, \cite{BD}, \cite{BP}, \cite{M}). Recently, the quenched invariance principle for  random conductance model with more general assumptions was shown  (\cite{ADS}, \cite{DNS}). They include not only the bond percolation cluster but also percolation clusters in models with long-range correlations and random conductance models defined by level sets of the Gaussian free field. In these papers, the proofs are based on an analysis on (general) weighted graph. They also use geometric conditions such as a relative isoperimetric inequality and a volume regularity. 
	
	In these results, they consider discrete settings. When we consider the continuum settings, the diffusion process is one of the fundamental objects. For the diffusion process, there are many homogenization results. In \cite{R} and \cite{R2}, annealed invariance principles were shown. (precisely, slightly stronger results were proved but not quenched results.)  For quenched results, there are \cite{CD}, \cite{FK}, \cite{FK2},\cite{FK3},\cite{O} and \cite{SZ}.   
	
	The homogenization problem of the continuum percolation cluster is considered by some researchers. Tanemura \cite{T} and Osada (\cite{O2}, \cite{OS}) proved the annealed invariance principle for the reflecting Brownian motion in the continuum percolation cluster. However, to our knowledge, there are no quenched results.

	\subsection{Method}
	A basic but powerful method to prove a quenched invariance principle is \textit{harmonic embedding}. A key ingredient of this approach is the \textit{corrector}, a random function, $ \chi\colon \mathbb{R}^d\times \Omega \to \mathbb{R}^d $ which is the solution of Poisson type equation
\begin{align}\label{eqn:har}
\mathcal{L}^\omega\chi^k(x,\omega) = \mathcal{L}^\omega\Pi^k(x), \; (k=1,\dots , d)
\end{align}
where $ \mathcal{L}^\omega $ is the corresponding generator of the diffusion $\{X_t^\omega\}_t$ and $ \Pi^k(x)=x_k $ is the projection to the $ k $-th coordinate. Then $ y^k = \Pi^k-\chi^k $ is a harmonic function. This implies that
\[ M_t^\omega = X_t^\omega-\chi(X_t^\omega,\omega)  \]
is a martingale and a quenched invariance principle for the martingale part $ M $ can be easily shown by standard arguments. In order to obtain a quenched invariance principle for the process $ X $, it suffices to show that for any $ T>0 $ and $ \HP $-a.s. $ \omega $
\[ \lim_{\epsilon \to 0 }\sup_{0 \leq t \leq T} \epsilon|\chi (X_{t/\epsilon^2}^\omega,\omega)| = 0 \quad \text{in $P_0^\omega$- probability}, \]
which can be deduced from the $ L^\infty  $-sublinearity of the corrector:
\begin{align}\label{for:slt}
	\lim _{\epsilon \to 0}\sup _{x \in W'_R}|\chi _\epsilon (x,\omega)|=0,\quad \HP \text{-a.s.},
\end{align}
where $\chi_\epsilon(x,\omega) = \epsilon\chi(x/\epsilon, \omega)$. The non-degeneracy of the covariance matrix follows from the $ L^\infty $-sublinearity and the ergodic theorem. So, important things are the followings:
\begin{enumerate}
	\item How to construct the corrector.
	\item How to prove the $ L^\infty $-sublinearity.
\end{enumerate}

One way to construct the corrector is to decompose the space of random function into a space of ``potential"  and its ``orthogonal'' space. 

Although the $ L^\infty $-sublinearity is difficult to show, a $ L^p $-sublinearity, which is the equation obtained from (\ref{for:slt}) replacing the $ L^\infty $-norm by the $ L^p $-norm, is easily shown by the ergodic theorem. In \cite{CD} and \cite{DNS}, it is mentioned that the $L^\infty $-sublinearity follows from the $ L^p $-sublinearity and a maximal inequality $(\ref{ine:max})$. They also suggested in these papers that Moser's iteration scheme is useful to obtain the maximal inequality. Since this method is analytic and robust, we can use this for the reflecting diffusion case.
	
	Throughout the proof, we will mainly follow the argument that appeared in Chiarini and Deuschel \cite{CD}. However, different from their paper, we need to consider the boundary effect. For example, we have to solve the equation (\ref{eqn:har}) with the Neumann boundary condition. Due to this boundary condition, the proof of the maximal inequality $(\ref{ine:max})$ is more complicated. To consider this boundary condition, we have to analyze the space $\CF^\omega$. Unlike the space $H_0^1(W)$ (the closure of 
	$ C_c^\infty(W) $ 
	with respect to a Sobolev norm), the trace of functions that belong to $\CF^\omega$ doesn't vanish in general. 
	Thus, we have to control the boundary effect to obtain the maximal inequality. To overcome this difficulty, we kill the process $X_t^\omega$ at the boundary $\partial W'(\omega,x,R) - \partial W'(\omega)$ and consider local inequalities in Section \ref{sec:maximal}.  

\section{Harmonic embedding}
For $x \in \mathbb{R}^d$ and $F\colon \Omega \to \mathbb{R}^d$, define $T_xF$ by $(T_xF)(\omega)=F(\tau _x \omega)$. Set \\ $\hl =\{F\colon \OH \to \mathbb{R}^d\mid \HE[ \langle aF,  F\rangle ]<\infty  \}$. We endow it with inner product $\langle aF, G \rangle$ for $F,G \in \hl$. 

We define the $i$-th derivative $D_i$ by
\[D_i U = \text{$L^2$-}\lim _{h \to 0}\frac{T_{he_i}U-U}{h},\]
where $e_i$ is the $i$-th coordinate vector and $U \in L^2(\OH ,\HP)$. 

Denote the domain of $D_i$ by $\mathcal{D}(D_i) = \{U \in L^2(\OH, \HP) \mid D_iU \text{ exists}  \}$. Set 
\[\C = \left \{ \int _{\mathbb{R}^d} f(\tau _x\omega)\phi (x)dx \setmid f \in  L^\infty (\Omega , \HP),  \phi \in C_c^\infty (\mathbb{R}^d) \right \}.\]
We can easily show that if $\displaystyle v(\omega)=\int_{\mathbb{R}^d} f(\tau _x\omega)\phi (x)dx \in \C$, then the $i$-th derivative of $ v $ exists and it is given by
	\begin{align}
		D_iv(\omega ) = -\int _{\mathbb{R}^d} f(\tau _x\omega ) \partial _i\phi (x)dx.
	\end{align}
Moreover, if $v \in \C$, then $v \in \bigcap _{i=1}^d \mathcal{D}(D_i)$.

\begin{Def} We define the gradient and the space of potentials as follows:
	\begin{enumerate}
	\item For $v \in \C$, set $D v=(D_1v, \dots , D_dv) \in \hl$.
	\item Define the subspace $ \hl_\text{pot} $ of $ \hl $ by the closure in $ \hl $ of the set $ \{Dv  \mid v \in \mathcal{C} \} $.
	\end{enumerate}
	\end{Def}
	
	\begin{rem}\label{rem:pot}
		\begin{enumerate}
			\item In the settings of \cite{CD}, the mean of a potential is zero and the authors of \cite{CD} used this property to prove the sublinearity of the corrector. Remind that the space $\mathcal{C}$ and the space of potentials can be generalized on $(\Omega, \mathbb{P})$. Then, for all generalized potential $\tilde{U}$, we have that $\E[\tilde{U}]=0$. However, the equality
			$\HE[U] = 0$
			does not hold for a general potential $ U \in \hl_\text{pot}$.
			\item For $ v \in \C $, set $v(x,\omega) = v(\tau_x\omega)  $. Then we have $ Dv(x,\omega) = \nabla v(x,\omega) $.
		\end{enumerate}
	\end{rem}
	\pagebreak
	\begin{lem}\label{lem:pot}
		Let $U \in \hl _{\text{pot}}$. Then for all $\eta \in C_c^\infty (\mathbb{R}^d)$ and $i,j=1,\dots, d$, we have
		\begin{align}
			\int _{\mathbb{R}^d}U_i(\tau _x\omega)\partial _j \eta (x)dx =\int _{\mathbb{R}^d}U_j (\tau _x\omega)\partial _i \eta (x)dx
		\end{align}
		for $\HP$-a.e. $\omega \in \Omega$. 
	\end{lem}
		
	\begin{proof}
		We first consider the case $v \in \C$.
		Then $x \mapsto v(\tau _x\omega)$ is infinitely many times differentiable, $\HP$-a.s. Integrating by parts we get 
		\begin{align*}
			\int _{\mathbb{R}^d}D_iv(\tau _x \omega)\partial _j\eta(x) dx  &= -\int _{\mathbb{R}^d}v(\tau _x \omega)\partial _i\partial _j \eta (x)dx \\
			&= -\int  _{\mathbb{R}^d}v(\tau _x \omega)\partial _j\partial _i \eta (x)dx \\
			&= \int _{\mathbb{R}^d}D_jv(\tau _x \omega)\partial _i \eta (x)dx.
		\end{align*}  
		For general $U \in \hl _{\text{pot}}$ take approximations and use the fact that as $n \to \infty $, $\nabla v_n \to U$ in $\hl$ implies $D_i v_n(\tau _\cdot \omega)\to U_i(\tau _\cdot \omega)$ in $L_{\text{loc}}^1$ $\HP$-a.s.
	\end{proof}

	Let $\pi ^k$ be the unit vector in the $k$-th direction. Since $\pi ^k \in \hl$, for each $k=1,\dots, d$, there exist functions $U^k \in \hl _{\text{pot}}$ and $R^k \in (\hl _{\text{pot}})^\perp$ such that $\pi ^k=U^k+R^k$.	Then we build the corrector starting from the functions $U^k \in \hl _{\text{pot}}$. For $k=1,\dots, d$, we define the corrector to be the function $\chi ^k \colon \mathbb{R}^d\times \Omega \to \mathbb{R}$ such that
	\[\chi ^k(x,\omega)=\sum _{j=1}^d \int _0^1 x_jU_j^k (\tau _{tx}\omega)dt.\]

	\begin{prop}[weak differentiability] For $k=1,\dots d$, the function $x \mapsto \chi ^k(x,\omega)$ is in $L_\text{loc}^1(\mathbb{R}^d)$, weakly differentiable $\P$-a.s. and $\partial _i \chi ^k(x,\omega)=U_i^k(\tau _x\omega)$ for $i = 1,\dots d$.
		\end{prop}
		\begin{proof}
		Let $\eta \in C_c^\infty (\mathbb{R}^d)$. We calculate
		\[\int _{\mathbb{R}^d}\chi ^k (x,\omega)\partial _i \eta (x)dx = \int _{\mathbb{R}^d}\left(\sum _{j=1}^d \int _0^1 x_jU_j^k(\tau _{tx}\omega)dt \right)\partial _i \eta (x)dx .\]
		By changing the order of integration (this can be done since $U_j^k (\tau _{tx}\omega) \in L_\text{loc}^1$) and applying the change of variables $y=tx$, we get
		\begin{align}\label{eqn:InPropWeakDel1}
			\int _{\mathbb{R}^d}\left(\sum _{j=1}^d \int _0^1 x_jU_j^k(\tau _{tx}\omega)dt \right)\partial _i \eta (x)dx = \int _0^1\sum _{j=1}^d \int _{\mathbb{R}^d} \frac{y_j}{t^{d+1}}U_j^k(\tau _y\omega)\partial _i \eta \left (\frac{y}{t} \right)dy dt 	
		\end{align}
		
		Since for $j\neq i$,
		\begin{align}\label{eqn:InPropWeakDel2}
			\frac{y_j}{t^{d+1}}\partial _i\eta \left(\frac{y}{t} \right)= \partial _i\left (\frac{y_j}{t^d}\eta \left(\frac{y}{t} \right) \right ),
		\end{align}
		Lemma \ref{lem:pot} and $(\ref{eqn:InPropWeakDel2})$ give
		\begin{align}\label{eqn:InPropWeakDel3}
			\int _{\mathbb{R}^d}\chi ^k (x,\omega)\partial _i \eta (x)dx =\int _{\mathbb{R}^d}U_i^k(\tau_y \omega)\int _0^1 \sum _{j\neq i}\partial _j\left (\frac{y_j}{t^d}\eta \left (\frac{y}{t} \right) \right) + \frac{y_i}{t^{d+1}}\partial _i \eta \left ( \frac{y}{t} \right) dtdy .
		\end{align}
		Finally, observe that for $y \neq 0$, 
		\begin{align}\label{eqn:InPropWeakDel4}
			\int _0^1 \sum _{j\neq i}\partial _j\left (\frac{y_j}{t^d}\eta \left (\frac{y}{t} \right) \right) + \frac{y_i}{t^{d+1}}\partial _i \eta \left ( \frac{y}{t} \right) dt= - \int _0^1 \frac{1}{t^{d-1}}\frac{d}{dt}\left (\eta \left (\frac{y}{t} \right ) \right)dt=-\eta (y)	
		\end{align}
		Combining $(\ref{eqn:InPropWeakDel3})$ and $(\ref{eqn:InPropWeakDel4})$, we have
		\[\int _{\mathbb{R}^d}\chi^k(x,\omega)\partial _i \eta (x)dx = - \int _{\mathbb{R}^d}U_i^k(\tau _x \omega)\eta (x)dx.\]
		Hence we have that for each $\eta \in C_c^\infty(\mathbb{R}^d)$, there exists a $\HP$-null set $N_\eta$ such that the above equality holds  for $\omega \in \OH \setminus N_\eta$. To obtain the desired result, we have to remove the dependency of $\eta$. Because $C_c^\infty (\mathbb{R}^d)$ is separable with respect to the supremum norm, we can remove this ambiguity considering a countable dense subset $\{\eta _n\}_n \subset C_c^\infty (\mathbb{R}^d)$ and a null set $N = \bigcup_n N_{\eta_n}$. 
		\end{proof}

		\begin{Def}
			We say that $u \in \CF_\text{loc}^\omega$ if for any relatively compact set $G \subset  \mathbb{R}^d$, there exists a function $u_G \in \CF^\omega$ such that $u=u_G$ a.e. on $G\cap W'(\omega)$.
			\end{Def}
			
		\begin{prop}
			For $k = 1, \dots d$, the corrector $\chi ^k (\cdot , \omega) \in \CF_\text{loc}^\omega $ for $\HP$-a.e. $\omega \in \Omega$.
			
		\end{prop}
			
		\begin{proof}
			By construction, there exists $\{f_n\}_n \subset \C$ such that $\nabla f_n \to U^k$, $n \to \infty$ in $\hl$. This implies that for any ball $B \subset \mathbb{R}^d$,
			\begin{align*}
			&\HE \left [\int _B \langle a(\tau_x\omega) \bigl(\nabla f_n (\tau _x \omega)-\nabla \chi (x, \omega)\bigr), \nabla f_n (\tau _x \omega)-\nabla \chi (x, \omega) \rangle dx \right ]  \\
			=&\int _B \HE [ \langle a(\tau_x\omega) \bigl(\nabla f_n (\tau _x \omega)-U^k (\tau_x\omega)\bigr), \nabla f_n (\tau _x \omega)-U^k (\tau_x\omega) \rangle ]dx \\
			\leq&\Lambda |B|\E[|\nabla f_n -U^k|^2]\P(0\in W')^{-1} \to 0  
			\end{align*} 
			Observe that $g_n(x,\omega)=f_n(\tau _x\omega)-f_n(\omega)$ belongs to $C_c^\infty (\mathbb{R}^d)$ and satisfies 
			\[g_n (x,\omega)=\sum _{i=1}^d \int _0^1 x_i \partial_i f_n(\tau _{tx} \omega)dt .\] 
			It follows that $g_n \to \chi ^k$ on $B$ with respect to $\norm{\cdot}_{L^2(B)}+\CE(\cdot,\cdot)$. This implies that $\chi ^k (\cdot, \omega)\in \CF_\text{loc}^\omega$ $\HP$-a.s.
		\end{proof}

		To obtain a martingale decomposition, we introduce a weak notion of the harmonicity.

		\begin{Def}
			We say that $u \in \CF_\text{loc}^\omega$ is $\CE^\omega$-harmonic if it satisfies 
			\[\CE^\omega (u,\phi)=0 \]
			for all $\phi \in C_c^\infty (\overline{W'(\omega)})$ .
			\end{Def}
			
			Set $y^k(x,\omega)=x_k-\chi ^k(x,\omega)$.

		\begin{prop}\label{prop:har}
			For $k=1,\dots, d$, $x \mapsto y^k(x,\omega)$ are $\CE^\omega $-harmonic $\P$-a.s.
		\end{prop}
		\begin{proof}
			It is enough to show that for any function $\phi \in C_c^\infty(\mathbb{R}^d)$,
			\begin{align}\label{equ:Eham}
				\int_{W'}\langle a\nabla y^k,\nabla \phi \rangle \;dx = 0
			\end{align}
			since a function in $C_c^\infty(\overline{W'(\omega)})$ is embedded in $C_c^\infty(\mathbb{R}^d)$.

			We introduce a subset of $ \C $. Define $ \C_0 $ by 
			\[ \C_0 = \left\{v \in \C \setmid v(\omega)=0 \text{ if } 0 \notin W'(\omega)   \right\} .\]
			Since 
			\[ L(\tau_x\omega) = \bigcup_{z\in \omega}B(z-x,\rho) =L(\omega)+x, \]
			we have 
			\[W'(\tau_x\omega) = W'(\omega)+x \] 
			and hence the condition $ 0 \notin W'(\tau_x\omega) $ is equivalent to $ -x \notin W'(\omega) $. 
			To prove (\ref{equ:Eham}), for all $f \in \mathcal{C}_0$, using Fubini's theorem we have 
			\begin{align*}
				&\HE\biggl[f(\omega)\int_{W'}\langle a(\tau_x\omega)\nabla y^k(x,\omega),\nabla \phi(x) \rangle \;dx \biggr] \\
				&= \HE\biggl[f(\omega)\sum_{i,j}\int_{W'} a_{ij}(\tau_x\omega)\partial_iy^k(x,\omega)\partial_j\phi(x)\;dx \biggr] \\
				&= \sum_{i,j}\int_{\mathbb{R}^d} \partial_j\phi(x)\HE [f(\omega)a_{ij}(\tau_x\omega)\partial_iy^k(x,\omega)\1{x\in  {W'}}]\;dx.
			\end{align*}
				Next using the shift invariance of $\P$, we have
			\begin{align*}
				&\sum_{i,j}\int_{\mathbb{R}^d} \partial_j\phi(x)\HE [f(\omega)a_{ij}(\tau_x\omega)\partial_iy^k(x,\omega)\1{x\in  W'}]\;dx \\
				&= \sum_{i,j}\int_{\mathbb{R}^d} \partial_j\phi(x) \E [f(\omega)a_{ij}(\tau_x\omega)\partial_iy^k(x,\omega)\1{x\in W'} \1{0\in W'}]\P(0\in W')^{-1}\;dx \\
				&= \sum_{i,j}\int_{\mathbb{R}^d} \partial_j\phi(x) \E [f(\tau_{-x}\omega)a_{ij}(\omega)\partial_iy^k(0,\omega)\1{0\in W'} \1{-x\in W'}]\P(0\in W')^{-1}\;dx\\
				&= \sum_{i,j}\HE\biggl[a_{ij}(\omega)\partial_iy^k(0,\omega) \int_{W'}f(\tau_{-x}\omega)\1{-x \in W'}\partial_j\phi(x)\;dx\biggr].
			\end{align*}
			Since $f \in \C_0$, the trace of a function $ x \mapsto f(\tau_x\omega)$ to $ \partial(-W'(\omega)) $ is identically zero. Hence by the Gauss-Green formula, we calculate
			\begin{align*}
				&\sum _{i,j=1}^d \HE \biggl [a_{ij}(\omega) \partial _i y^k (0,\omega) \int_{-W'} f(\tau_x\omega) \partial _j \phi (-x)dx  \biggr ] \\
				&=  -\sum _{i,j=1}^d \HE \biggl [a_{ij}(\omega)  \partial _i y^k (0,\omega) \int_{-W'} \partial _j f(\tau_x\omega) \phi (-x)dx  \biggr ] \\
				&=  -\sum _{i,j=1}^d \HE \biggl [a_{ij}(\omega)  \partial _j y^k (0,\omega) \int_{W'} \partial _i f(\tau_{-x}\omega) \phi (x)dx  \biggr ]. \\ 
			\end{align*} 
			Again, by  the shift invariance, we find that the last line equals to
			\begin{align*}
				-\sum _{i,j=1}^d \int _{\mathbb{R}^d}\phi (x) \HE [a_{ij}(\omega) \partial _i y^k (0,\omega)  \partial _i f(\tau _{-x}\omega)\1{-x \in W'} ] dx
			\end{align*}
			Because $f$ belongs to $\C_0$, the function $ x \mapsto \partial_i f(\tau_x\omega) $ vanishes on $ \mathbb{R}^d\setminus (\overline{-W'(\omega)}) $. Hence we calculate 
			\begin{align*}
				&= -\sum _{i,j=1}^d \int _{\mathbb{R}^d}\phi (x) \HE [a_{ij}(\omega) \partial _i y^k (0,\omega)  \partial _i f(\tau _{-x}\omega)\1{-x \in W'} ] dx \\
				&= -\sum _{i,j=1}^d \int _{\mathbb{R}^d}\phi (x) \HE [ a_{ij}(\omega)\partial _i y^k (0,\omega)  \partial _i f(\tau _{-x}\omega) ] dx \\
				&= \int _{\mathbb{R}^d} \phi(x) \HE [\langle a(\omega)\nabla y^k(0,\omega), \nabla f(\tau_{-x}\omega)\rangle] dx.
			\end{align*}
			Since $ \tau_{-x}\C \subset \C $ and $\nabla y^k \in (\hl _{\text{pot}})^\perp$, we get
			\begin{align*}
				\int _{\mathbb{R}^d} \phi(x) \HE [\langle a\nabla y^k(0,\omega),\nabla f(\tau_{-x}\omega)\rangle ] dx =0.
			\end{align*}
			Since $\C_0 \subset L^p (\Omega ,\HP)$ for all $p\geq 1$ densely, it follows that 
			\[\int _{W'}\langle a\nabla y^k (x,\omega), \nabla \phi (x) \rangle dx =0, \;\; \HP\text{-a.s.} \]
			This ends the proof. 
		\end{proof}

		If a function  $ f $ is harmonic, the process $ M_t = f(X_t) $ is a martingale from a standard result of Markov processes. To deduce the same result for $ \CE^\omega  $-harmonic functions, we will use the following theorem due to Fukushima, Nakao and Takeda \cite[Theorem 3.1]{FNT}. We refer the terminology of Dirichlet forms to \cite{FOT}.

		Let $W$ be a Lipschitz domain on $\mathbb{R}^d$. Let $(Y_t,Q_x)$ be a Hunt process associated with strongly local regular Dirichlet form $(\CE  , \CF)$ on $ L^2(W)$. Fix a point $x_0 \in W$ and consider the following conditions for the process $(Y_t, Q_x)$ and for a function $u$:
	\begin{enumerate}
			\item The transition function $p_t$ of $X_t$ satisfies $p_t(x_0,A)=0$ for any $t$ if $\cp{A}=0$.
			\item The function $u$ belongs to $\CF_\text{loc}$. Moreover, the function $u$ is continuous and $\CE$-harmonic.
			\item The energy measure $\nu  _{\langle u\rangle}$ of $u$ is absolutely continuous with respect to the Lebesgue measure on $W$ and the density function $f$ satisfies
			\[E_{x_0}\biggl [\int _0^t f(X_s)dx\biggr ]<\infty, \;\; t>0 .\]
	\end{enumerate}
\begin{lem}(\cite[Theorem 3.1]{FNT})\label{lem:mart}
	Assume that the above conditions hold. Then the additive functional $M_t = u(X_t)-u(X_0)$ is a $P_{x_0}$-square integrable martingale with
\begin{align}
\langle M \rangle _t = \int _0^t f(X_s)ds, \;\; t>0, \;\; P_{x_0}\text{-a.s.}  
\end{align}
\end{lem}

\begin{cor}\label{cor:mart}
	The process $y (X_t, \omega)$ is a martingale with quadratic variation
	\[\langle y^k(X_\cdot,\omega), y^\ell(X_\cdot, \omega) \rangle_t = \int _0^t \sum _{i,j=1}^da_{ij}(\tau_x\omega)\partial _iy^k(X_t,\omega)\partial _jy^\ell(X_t,\omega)ds, \; (k,\ell = 1,\dots , d)  .\] 
\end{cor}

\begin{proof}
	To apply Lemma \ref{lem:mart} for $y^k$, we need to check the assumptions of Lemma $\ref{lem:mart}$.
	
	Since the reflecting diffusion $\{X_t\}_t$ has the transition density (say $ q_t(x,y) $), we have
	$p_t(x_0,A) =\int _A q_t(x_0,y)dy = 0$ 
	if $\cp{A}=0$. Hence assumption (1) of Lemma \ref{lem:mart} holds.
	
	By Proposition \ref{prop:har}, for $ \HP $-a.s. \hspace{-1em} $ \omega $, the function $ x \mapsto y^k(x,\omega) $ satisfies a solution of the equation
	\[\int _W \langle a(\tau_x\omega)\nabla y^k(x, \omega), \nabla \phi(x)\rangle dx = 0 , \quad \phi \in C_c^\infty (W)   .\]  
	Thus, assumption (2) follows from classical results in elliptic partial differential equations (see \cite[Section 9]{GT}). We have
	\begin{align*}
	&2\CE^\omega (y^kv,y^k)-\CE^\omega ((y^k)^2,v)\\
	&= 2\int_{W'} v\langle a(\tau_x\omega)\nabla y^k, \nabla y^k\rangle  dx + 2\int_{W'} y^k\langle a(\tau_x\omega)\nabla v, \nabla y^k\rangle dx - 2\int_{W'} y^k\langle a(\tau_x\omega)\nabla y^k, \nabla v \rangle dx \\
	&= \int_{W'} 2\langle a(\tau_x\omega)\nabla y^k, \nabla y^k\rangle v dx
	\end{align*}
	for all $ v \in C_c^\infty(W) $.
	Since the function $ x \mapsto y^k(x,\omega) $ is weakly differentiable, the energy measure $ \nu_{\langle y^k \rangle} $ is given by 
	\[ \nu_{\langle y^k \rangle}(dx) = 2\langle a(\tau_x\omega)\nabla y^k(x), \nabla y^k(x)\rangle dx \]
	and absolutely continuous with respect to the Lebesgue measure. Moreover, its density $ f $ is written by 
	\[ f(x) = 2\langle a(\tau_x\omega)\nabla y^k(x), \nabla y^k(x)\rangle. \]
	Next using the stationarity of the environment process $ \{\tau_{X_t^\omega}\}_t $ under $ \HP $ (see \cite{O}), by using Fubini's Theorem and the translation invariance, we compute 
	\begin{align*}
	\HE \left [ E_0^\omega \left[\int_0^t f(X_s^\omega)ds \right] \right ]
	= \int_0^t\HE[E_0^\omega[f(X_s^\omega)]]ds  
	= \int_0^t\HE [E_0^\omega[f(X_0^\omega)]]ds 
	= t\HE[E_0^\omega[f(0)]]< \infty. 
	\end{align*}
	Therefore, assumption (3) is satisfied. 
	\end{proof}

\section{Maximal inequality}\label{sec:maximal}
In this section, we consider deterministic settings.
Let $ W \subset \mathbb{R}^d $ be an unbounded Lipschitz domain containing the origin. We denote the connected component of $ B(0,R)\cap W $ containing the origin by $ W_R $. We assume that there exists positive constants $C_V$ and $R_V$ such that 
\begin{align}\label{ine:vol}
|W_R|\geq C_V R^d
\end{align}
for $R \geq R_V$. We further assume that there exist $\theta \in (0,1)$ and $c_H>0$ such that  
\begin{gather}\label{cdn:InThmSobMIsoperi}
	C_\text{IL} := \inf \left\{\frac{\CH_{d-1}(W \cap \partial O)}{|W \cap O |^\frac{d-1}{d}} \middle | \begin{split}
		 &\text{ $O\subset B(0,R)$ is connected open,} \\
		 &R\geq R_I,\;|O| \geq R^\theta  \end{split} \right\} > 0
\end{gather}
for some $R_I>0$ and 
\begin{gather}\label{cdn:IPforShort}
	C_ \text{IS}:=\inf\left\{\frac{\mathcal{H}_{d-1}(W\cap \partial O)}{|W\cap O|^\frac{d-1}{d}}\middle | \begin{split} &O\text{ is bounded open,}\\ &\mathcal{H}_{d-1}(W\cap \partial O) < c_H
	\end{split}\right\} > 0
\end{gather}
hold.


Set $\zeta = \frac{1-\theta}{1-\frac{\theta}{d}}$. For a bounded open subset, the following weak isoperimetric inequality holds.

\begin{lem}\label{lem:WI}
    There exists positive constant $C_I$ such that 
    \begin{align}\label{ine:wiso}
        \frac{\mathcal{H}_{d-1}(W\cap \partial O)}{|O|^\frac{d-\zeta}{d}}\geq \frac{C_I}{R^{1-\zeta}}
    \end{align}
    holds for open subset $O \subset W(0,R)$ and $R \geq R_I$.
\end{lem}
\begin{proof}
 The general idea of the proof goes back to \cite[Lemma, 3.3]{CD}. First we consider a good and connected open subset $O$. When $O$ satisfies $|O|\geq R^\theta$ or $\mathcal{H}_{d-1}(W\cap \partial O)<c_H$, from $(\ref{cdn:InThmSobMIsoperi})$ and $(\ref{cdn:IPforShort})$, we have 
 \begin{align}\label{ine:inWI1}
	\frac{\mathcal{H}_{d-1}(W\cap\partial O)}{|O|^\frac{d-\zeta}{d}}=\frac{\mathcal{H}_{d-1}(W\cap\partial O)}{|O|^\frac{d-1}{d}|O|^\frac{1-\zeta}{d}}\geq \frac{C_\text{IL}\wedge C_\text{IS}}{|O|^\frac{1-\zeta}{d}}\geq\frac{C_\text{IL}\wedge C_\text{IS}}{|B(0,R)|^\frac{1-\zeta}{d}}\geq\frac{C}{R^{1-\zeta}}.
 \end{align}
 Next we consider the case where $|O|< R^\theta$ and $\mathcal{H}_{d-1}(W\cap \partial O) \geq c_H$. Using the relation $\theta \cdot \frac{d-\zeta}{d} = 1-\zeta$, we estimate 
 \begin{align}\label{ine:inWI2}
	\frac{\mathcal{H}_{d-1}(W\cap\partial O)}{|O|^\frac{d-\zeta}{d}}\geq \frac{c_H}{|O|^\frac{d-\zeta}{d}}\geq\frac{c_H}{R^{\theta \cdot \frac{d-\zeta}{d}}} = \frac{c_H}{R^{1-\zeta}}.
 \end{align}
 Combining $(\ref{ine:inWI1})$ and $(\ref{ine:inWI2})$, 
 we get $(\ref{ine:wiso})$ when $O$ is connected. For a general open subset $O$, it follows from the fact that $(a_1+a_2)/(b_1+b_2)^\beta \geq a_1/b_1^\beta+a_2/b_2^\beta \geq M$ for $\beta\in (0,1)$ and positive numbers $a_1,a_2,b_1,b_2,M$  with $a_1/b_1^\beta\geq M$ and $a_2/b_2^\beta\geq M$. 
\end{proof}

Let $E\subset W$ be a bounded set. For a function $u\colon E \to \mathbb{R}$, $\alpha \geq 1$ we denote 

\begin{align*}
	\norm{u}_{L^\alpha(E)} := \left( \int_E |u(x)|^\alpha dx \right)^\frac{1}{\alpha},\; \norm{u}_{E,\alpha} := \left(\frac{1}{|E|}\int_E |u(x)|^\alpha dx\right)^\frac{1}{\alpha}.
\end{align*}

Next we show a Sobolev type inequality. To do this, we need to consider functions whose trace vanishes locally. For $R>0$, set
\begin{align*}
	\mathcal{C}_R = \{u \in C_c^\infty(\overline{W_R}) \mid u = 0 \text{ on $\partial W_R - \partial W$} \}.
\end{align*}
Let $ a\colon \mathbb{R}^d \to \mathbb{R}^{d\times d} $ be a positive-definite symmetric matrix. Define the bilinear form $ \CE $ on $L^2(W, dx)$ by
\begin{align}
\CE (u,v)=\int_W \langle a \nabla u, \nabla v \rangle dx, \quad u,v \in L^2(W,dx)\cap C^\infty(W).
\end{align}
\begin{asm}\label{asm:Dir} 
	\begin{enumerate}
		\item there exist positive measurable functions $ \lambda, \Lambda $ such that for almost all $x, \xi \in \mathbb{R}^d$,
		\[ \lambda (x)|\xi|^2 \leq \langle a(x)\xi, \xi \rangle \le \Lambda (x)|\xi|^2 .\]
		\item there exist $ p,q \in [1,\infty] $ satisfying $ 1/p+1/q < 2\zeta/d $ such that for almost all $x$,
		\begin{align*}
		\limsup_{R \to \infty} \frac{1}{|W_R|}\int_{W_R}(\Lambda^p+\lambda^{-q})dx < \infty.
		\end{align*}
	\end{enumerate}
\end{asm}
We remark that the relation between $ p $ and $ q $ is more restricted than Chiarini and Deuschel \cite{CD}. This is because of the boundary effect.
Let $\CF_R$ be the closure of $\mathcal{C}_R$ with respect to $\norm{\cdot}_{L^2(\C_R)}+\CE(\cdot,\cdot)$. 
Then if the Dirichlet form $(\CE,\CF_R)$ is regular, the associate diffusion is absorbed at $\partial W_R - \partial W$ and has a reflection at $\partial W_R\cap \partial W$.
We note that the domain $\mathcal{C}_R$ is larger than the Sobolev space with zero boundary condition $H_0^1(W_R)$. Hence, it is not obvious whether the Sobolev inequality holds.  The key property is the weak isoperimetric inequality 
$(\ref{ine:wiso})$.

\begin{prop}\label{pro:Sob1}
    Let $R\geq R_I \vee R_V$.
    Then we have 
    \begin{align}\label{ine:sobM}
        \norm{u}_{L^\frac{d}{d-\zeta}(W_R)} \leq C_S^{-1} |W_R|^\frac{1-\zeta}{d} \norm{\nabla u}_{L^1(W_R)}
    \end{align}
    for $u \in \CF_R$, where $C_S=C_IC_V^\frac{1-\zeta}{d}$.
\end{prop}
\begin{proof}
    We first show $(\ref{ine:sobM})$ for non-negative $u \in \mathcal{C}_R$.
    Let $v \colon \mathbb{R}^d \to \mathbb{R}$ be a zero-extension of $u$. That is, $v = u$ in $\text{supp}\; u$ and $v = 0$ outside the $\text{supp}\; u$. 
    Then by the coarea formula (see \cite[Theorem 3.10]{E}), we have
    \begin{align*}
        \int_{W_R} |\nabla v| \;dx = \int_{-\infty}^\infty\CH_{d-1}(W_R \cap v^{-1}(\{t\}))\; dt.
    \end{align*}
    Now, since $v$ is non-negative, we have $v^{-1}(\{t\}) = \emptyset$ for $t < 0$. Moreover, $u = v$ on $W_R$. Thus, we have
    \begin{align}\label{ine:InsobMCor}
        \int_{W_R} |\nabla u| \;dx = \int_0^\infty\CH_{d-1}(W_R \cap \partial\{x \in W_R \mid u(x) > t\})\; dt.
    \end{align}\label{eqn:coarea0}
    Since $u$ is continuous and $\min u = 0$, we have $W_R \cap \{x \in W_R \mid u(x) = t\} = W_R \cap \partial \{x \in W_R \mid u(x) > t\}$ for $t > 0$. Since $u = 0$ on $\partial W_R - \partial W$, we also have $W_R \cap \partial \{x \in W_R \mid u(x) > t\} = W \cap \partial \{x \in W_R \mid u(x) > t\}$. Combining these and $(\ref{ine:InsobMCor})$, we obtain
    \begin{align}\label{eqn:coarea}
        \int_{W_R}|\nabla u|\; dx = \int_0^\infty\CH_{d-1}(W \cap \partial \{x \in W_R \mid u(x) > t\}) \; dt.
    \end{align}
    Set $U_t = \{x \in W_R \mid u(x) > t\}$. From  $(\ref{ine:wiso})$,  we estimate
    \begin{align}\label{ine:InsobMAfterCor}
        &\int_0^\infty\CH_{d-1}(W \cap \partial \{x \in W_R \mid u(x) > t\}) \; dt \nonumber \\
        &\geq \frac{C_I}{R^{1-\zeta}} \int_0^\infty |W \cap \{x \in W_R \mid u(x) > t\}|^\frac{d-\zeta}{d} \; dt \nonumber \\
        &= \frac{C_I}{R^{1-\zeta}}\int_0^\infty |\{x \in W_R \mid u(x) > t\}|^\frac{d-\zeta}{d} \; dt \nonumber\\
        &= \frac{C_I}{R^{1-\zeta}}\int_0^\infty |U_t|^\frac{d-\zeta}{d} \; dt \nonumber \\
        &= \frac{C_I}{R^{1-\zeta}}\int_0^\infty \norm{1_{U_t}}_{L^\frac{d}{d-\zeta}(W_R)} \; dt.
    \end{align}
    Using $(\ref{ine:vol})$, we have $1/R^{1-\zeta} \geq (C_V/|W_R|)^\frac{1-\zeta}{d}$. Combining this with $(\ref{eqn:coarea})$ and $(\ref{ine:InsobMAfterCor})$, we obtain
    \begin{align}\label{ine:InSob2nd}
        \int_{W_R}|\nabla u|\; dx  \geq \frac{C_S}{|W_R|^\frac{1-\zeta}{d}}\int_0^\infty \norm{1_{U_t}}_{L^\frac{d}{d-\zeta}(W_R)} \; dt.
    \end{align}
   Now let $r$ be positive number satisfying $(d/(d-\zeta))^{-1} + 1/r = 1$. Take $g \in L^r$ such that $g \geq 0$ and $\norm{g}_{L^r(W_R)}=1$. Then, by the H\"{o}lder inequality,
   \begin{align*}
       \int_0^\infty \norm{1_{U_t}}_{L^\frac{d}{d-1}(W_R)}dt \geq \int_0^\infty\norm{g 1_{U_t}}_{L^1(W_R)} \;dt 
       = \int_{W_R} g(x) \int_0^\infty 1_{U_t}(x)\;dtdx 
       = \norm{gu}_{L^1(W_R)}.
   \end{align*}
   Since $g$ is arbitrary, putting this into $(\ref{ine:InSob2nd})$, we obtain 
   \begin{align*}
       \norm{u}_{L^\frac{d}{d-\zeta}(W_R)} = \sup_{g \in L^r(W_R), \;  g \geq 0, \; \norm{g}_{L^r(W_R)} = 1} \norm{gu}_{L^1(W_R)} \leq C_S^{-1}|W_R|^\frac{1-\zeta}{d}\norm{\nabla  u}_{L^1(W_R)},
   \end{align*}
   hence $(\ref{ine:sobM})$ holds. For general $u \in \CF_R$, we can show the same bound by approximation.
\end{proof}
Set $p^* = \frac{p}{p-1}, \;\rho = \frac{2qd}{q(d-2\zeta)+d}$, and $q^\#=\frac{2q}{q+1}$.  

\begin{prop}\label{prop:Sob}
	Let $R\geq R_I\vee R_V$. Let $p,q$ be positive numbers as in (2) of Assumption \ref{asm:Dir}. Then we have 
	\begin{align}\label{ine:Sob}
		\norm{u}_{L^\rho(W_R)} \leq C_\text{sob}|W_R|^\frac{1-\zeta}{d}\norm{\nabla u}_{L^{q^\#}(W_R)}
	\end{align} 
	for $u \in \CF_R$.
\end{prop}
\begin{proof}
    First observe that $\rho = \frac{dq^\#}{d-q^\#\zeta}$.
    Let $r$ be a real number satisfying $1/q^\# + 1/r = 1$. Set $\bar{d}= q^\#(d-\zeta)/(d-q^\#\zeta)$. Using a relation $\frac{dq^\#}{d-q^\#\zeta} = \frac{q^\#(d-\zeta)}{d-q^\#\zeta}\frac{d}{d-\zeta}=\bar{d}\cdot \frac{d}{d-\zeta}$ and  Proposition \ref{pro:Sob1}, we have
    \begin{align}\label{ine:usingSob1}
        &\left(\int_E |u|^{\frac{dq^\#}{d-q^\#\zeta}} dx \right)^{\frac{d-1}{d}} =\norm{u^{\bar{d}}}_{L^\frac{d}{d-\zeta}(W_R)} \nonumber\\
        &\leq C_S^{-1}|W_R|^\frac{1-\zeta}{d}\norm{\nabla u^{\bar{d}}}_{L^1(W_R)} 
		=\bar{d}C_S^{-1}|W_R|^\frac{1-\zeta}{d}\norm{u^{\bar{d}-1}\nabla u}_{L^1(W_R)}. 
    \end{align}
	Using the H\"{o}lder inequality with a pair $(r,q^\#)$ and a relation $(\bar{d}-1)r = \frac{dq^\#}{d-q^\#\zeta}$, we have
	\begin{align*}
		&\bar{d}C_S^{-1}|W_R|^\frac{1-\zeta}{d}\norm{u^{\bar{d}-1}\nabla u}_{L^1(W_R)} \leq \bar{d}C_S^{-1}|W_R|^\frac{1-\zeta}{d}\norm{u^{\bar{d}-1}}_{L^r(W_R)}\norm{\nabla u}_{L^{q^\#}(W_R)} \\
        &= \bar{d}C_S^{-1}|W_R|^\frac{1-\zeta}{d}\norm{\nabla u}_{L^{q^\#}(W_R)}\left(\int_{W_R} |u|^\frac{dq^\#}{d-q^\#\zeta} dx\right)^\frac{q^\#-1}{q^\#}.
	\end{align*}
    Inserting this into $(\ref{ine:usingSob1})$, we get 
	\begin{align*}
		\left(\int_{W_R} |u|^{\frac{dq^\#}{d-q^\#\zeta}} dx \right)^{\frac{d-1}{d}}
		\leq \bar{d}C_S^{-1}|W_R|^\frac{1-\zeta}{d}\norm{\nabla u}_{L^{q^\#}(W_R)}\left(\int_{W_R} |u|^\frac{dq^\#}{d-q^\#\zeta} dx\right)^\frac{q^\#-1}{q^\#}.
	\end{align*}
	Dividing both side by $\displaystyle \left(\int |u|^\frac{dq^\#}{d-q^\#} dx\right)^\frac{q^\#-1}{q^\#},$ we get the desired result since $(d-\zeta)/d - (q^\#-1)/q^\# = (d-q^\#\zeta)/(dq^\#)$.
\end{proof}

Thanks to Proposition \ref{prop:Sob}, we can prove the following inequality.

\begin{lem}\label{lem:locSob}
	Let $R\geq R_I \vee R_V$. Then for all $u \in \CF_R$
	\begin{align}
		\norm{u}_{L^\rho(W_R)}^2 \leq C_\text{sob}|W_R|^\frac{1-\zeta}{d}\norm{\lambda^{-1}}_{L^q(W_R)}\CE(u,u).
	\end{align}
\end{lem}
\begin{proof}
    By  $(\ref{ine:Sob})$, we have 
    \begin{align*}
        \norm{u}_{L^\rho(W_R)} \leq C_\text{sob}|W_R|^\frac{1-\zeta}{d}\norm{\nabla u}_{L^{q^\#}(W_R)}.
    \end{align*}
    By the H\"{o}lder inequality and (1) of Assumption \ref{asm:Dir}, the right hand side is estimated by
\begin{align*}
\norm{\nabla u}_{L^{q^\#}(W_R)} &= \left(\int_{W_R} |\nabla u|^\frac{2q}{q+1}\lambda^\frac{q}{q+1}\cdot \lambda^{-\frac{q}{q+1}} \; dx\right)^\frac{q+1}{2q} \\
& \leq \left(\int_{W_R} (|\nabla u|^2\lambda)^{\frac{q}{q+1}\cdot\frac{q+1}{q}}dx\right)^{\frac{q+1}{2q}\cdot \frac{q}{q+1}} \norm{\lambda^{-1}}_{L^q(W_R)} \\
&\leq\CE (u,u)^\frac{1}{2}\norm{\lambda^{-1}}_{L^q(W_R)}.
\end{align*}
Therefore, we get the desired result.
\end{proof}


By cutoff on $W_R$ we mean a function $\eta \in \mathcal{C}_R$ satisfying  $0 \leq \eta \leq 1$.

\begin{prop}
    Let $R\geq R_I \vee R_V$. Let $ \eta $ be a cutoff on $W_R$. Then there exists a constant $C>0$, depending only on the dimension $d \geq 2$, such that for all $u \in \CF_R $ 
    \begin{align}\label{ine:cut}
    \norm{\eta u}_{L^\rho(W_R)}^2 \leq 2C |W_R|^\frac{1-\zeta}{d}\norm{\lambda^{-1}}_{L^q(W_R)}(\CE_\eta (u,u)+\norm{\eta}_\infty^2\norm{u\Lambda^\frac{1}{2}}_{L^2(W_R)}^2),
    \end{align} 
    where we denote $\displaystyle \CE _\eta(u,u)= \int_W \langle a\nabla u, \nabla u\rangle \eta^2 \;dx $.
    \end{prop}
\begin{proof}
By Lemma \ref{lem:locSob}, we have
\[\norm{\eta u}_\rho^2 \leq C_\text{sob} |W_R|^\frac{1-\zeta}{d}\norm{\lambda^{-1}}_{L^q(W_R)}\CE(\eta u,\eta u).\]
Because $\langle a(\eta\nabla u + u\nabla \eta), (\eta\nabla u + u\nabla \eta) \rangle \leq 2(\eta^2\langle a\nabla u,\nabla u\rangle + u^2\langle a\nabla\eta, \nabla\eta\rangle)$, we estimate 
\begin{align*}
    \CE(\eta u,\eta u) &= \int_W \langle a\nabla (\eta u), \nabla (\eta u)\rangle  \;dx \\
    &= \int_W \langle  a(\eta\nabla u + u\nabla \eta), (\eta\nabla u + u\nabla \eta) \rangle  \;dx \\
    &\leq \int_W 2(\eta^2\langle a\nabla u,\nabla u\rangle + u^2\langle a\nabla\eta, \nabla\eta\rangle)\;dx \\
    &\leq 2\int_W \langle a\nabla u , \nabla u\rangle \eta^2 \;dx + 2 \int_W \norm{\nabla\eta}_\infty^2 u^2 \Lambda dx \\
    &= 2\CE_\eta(u,u)+2\norm{\eta}_\infty^2\norm{u\Lambda^\frac{1}{2}}_{L^2(W_R)}^2,
\end{align*}
which leads to the conclusion.
\end{proof}

Let $\CF$ be the closure of $C_c^\infty(\overline{W})$ with respect to $\norm{\cdot}_{L^2(W)}+\CE(\cdot,\cdot)$. We say that a function $u$ belongs to $\CF_\text{loc}$ if for all $R>0$ there exists $u_R \in \CF$ such that $u = u_R$ on $W_R$. Let $f\colon W\to\mathbb{R}$ be a function with essentially bounded derivatives. Consider the following equation: 
\begin{align}\label{eqn:Poisson}
	\CE(u,\phi) = -\int_W \langle a\nabla f, \nabla \phi \rangle \; dx. 
\end{align}
We say that $u \in \CF_\text{loc}$ is a solution of the  equation $(\ref{eqn:Poisson})$ if it
holds for all $\phi  \in C_c^\infty(\overline{W})$. 
We say that $u$ is a subsolution of the equation $(\ref{eqn:Poisson})$ if the equation $(\ref{eqn:Poisson})$ holds with $\leq$ for all $\phi  \in C_c^\infty(\overline{W})$. 
We also say that $u$ is a solution of the equation in $W_R$ if the equation $(\ref{eqn:Poisson})$ holds for $\phi \in \CF_R$.

\begin{prop}
    Let $R\geq R_I \vee R_V$ and $u \in \CF_\text{loc}$ be a subsolution of equation of $(\ref{eqn:Poisson})$ in $W_R$. Let $\eta \in C_c^\infty(W_R)$ be a cutoff. Then there exists a constant $C_1 > 0$ such that for all $\alpha \geq 1$, 
	\begin{align}
		\norm{\eta u}_{W_R,\alpha \rho}^{2\alpha} \leq& \alpha ^2 C_1\norm{\lambda^{-1}}_{W_R,q}\norm{\Lambda}_{W_R,p}|W_R|^\frac{2}{d} \\
		&\times (\norm{\nabla \eta}_\infty^2 \norm{u^+}_{W_R, 2\alpha p^*}^{2\alpha} + \norm{\nabla f}_\infty^2 \norm{u^+}_{W_R,2\alpha p^*}).  \nonumber
	\end{align}
\end{prop}
    \begin{proof}
    Similarly to \cite[Proposition 2.4]{CD}, we can show that
    \begin{align*}
    \norm{\eta(u^++\epsilon)^\alpha}_\rho^2 
    &\leq 2 C |W_R|^\frac{1-\zeta}{d}\norm{\lambda^{-1}}_{L^q(W_R)}\bigl[(\alpha^2 +1)\norm{(u^+ + \epsilon)^{2\alpha}\Lambda}_{L^1(W_R)}\norm{\nabla \eta}_{L^\infty(W_R)}^2 \\ 
    &+ \norm{\nabla f}_{L^\infty(W_R)}^2 \alpha^2 \norm{(u^++\epsilon)^{2\alpha-2}\Lambda}_{L^1(W_R)} \\
    &+\frac{\alpha^2}{2\alpha-1}\norm{\nabla \eta}_{L^\infty(W_R)}\norm{\nabla f}_{L^\infty(W_R)}\norm{u^{2\alpha-1}\Lambda}_{L^1(W_R)} \bigr].
    \end{align*}
    (The only difference is using $(\ref{ine:cut})$ instead of \cite[Proposiotion 2.3]{CD}.).
    Taking the limit as $\epsilon \to 0$ and using the H\"{o}lder inequality with $1/p+1/p^*=1$, we have
    \begin{align*}
			\norm{\eta(u^+)^\alpha}_\rho^2 
        &\leq 2 C |W_R|^\frac{1-\zeta}{d}\norm{\lambda^{-1}}_{L^q(W_R)}\bigl[(\alpha^2 +1)\norm{(u^+)^{2\alpha}}_{L^{p^*}(W_R)}\norm{\Lambda}_{L^p(W_R)}\norm{\nabla \eta}_{L^\infty(W_R)}^2 \\ 
        &+ \norm{\nabla f}_{L^\infty(W_R)}^2 \alpha^2 \norm{(u^+)^{2\alpha-2}\Lambda}_{L^{p^*}(W_R)}\norm{\Lambda}_{L^p(W_R)}\\
        &+\frac{\alpha^2}{2\alpha-1}\norm{\nabla \eta}_{L^\infty(W_R)}\norm{\nabla f}_{L^\infty(W_R)}\norm{u^{2\alpha-1}}_{L^{p^*}(W_R)}\norm{\Lambda}_{L^p(W_R)} \bigr].
    \end{align*}
    
    Averaging over $W_R$ and using the relation $2/\rho = 1/p^* + 1/p +1/q - (2\zeta)/d$ and $(1+\zeta)/d \leq 2/d$, we get 
    \begin{align*}
    \norm{\eta(u^+)^\alpha}_{W_R,\rho}^2 &\leq 2 C \norm{\lambda^{-1}}_{W_R,q}\norm{\Lambda}_{W_R,p} |W_R|^\frac{1+\zeta}{d}\bigl[ (\alpha ^2 +1)\norm{(u^+ )^{2\alpha}}_{W_R,p^*}\norm{\nabla \eta}_\infty^2 \\ 
    &+ \alpha ^2 \norm{\nabla f}_\infty^2\norm{(u^+)^{2\alpha-2}}_{W_R,p^*}+\frac{\alpha^2}{2\alpha-1}\norm{\nabla \eta}_\infty\norm{\nabla f}_\infty\norm{(u^+)^{2\alpha-1}}_{W_R,p^*} \bigr]\\
    &\leq 2 C \norm{\lambda^{-1}}_{W_R,q}\norm{\Lambda}_{W_R,p} |W_R|^\frac{2}{d}\bigl[ (\alpha ^2 +1)\norm{(u^+ )^{2\alpha}}_{W_R,p^*}\norm{\nabla \eta}_\infty^2 \\ 
    &+ \alpha ^2 \norm{\nabla f}_\infty^2\norm{(u^+)^{2\alpha-2}}_{W_R,p^*}+\frac{\alpha^2}{2\alpha-1}\norm{\nabla \eta}_\infty\norm{\nabla f}_\infty\norm{(u^+)^{2\alpha-1}}_{W_R,p^*} \bigr].
    \end{align*}
    By Jensen's inequality we have 
    \[\norm{u^+}_{W_R,(2\alpha-2)p^*} \leq \norm{u^+}_{W_R, 2\alpha p^*}, \quad\norm{u^+}_{W_R,(2\alpha-1)p^*} \leq \norm{u^+}_{W_R, 2\alpha p^*}  ,\]
    therefore we can rewrite and get 
    \begin{align*}
    \norm{\eta u^+}_{W_R,\rho}^2 
    &\leq 2 C |W_R|^\frac{2}{d}\bigl[ (\alpha ^2 +1)\norm{1_{W_R}u^+ }_{W_R,2 \alpha p^*}^{2\alpha}\norm{\nabla \eta}_\infty^2 \\ 
    &+ \alpha ^2 \norm{\nabla f}_\infty^2\norm{1_{W_R}u^+}_{W_R,2 \alpha p^*}^{2\alpha -2}+\frac{\alpha^2}{2\alpha-1}\norm{\nabla \eta}_\infty\norm{\nabla f}_\infty\norm{u^+}_{W_R,2 \alpha p^*}^{2\alpha-1} \bigr].
    \end{align*} 
    Absorbing the mixed product in the two squares, we obtain the desired result. 
    \end{proof}

We can prove the following inequality as in \cite[Corollary 2.1]{CD}. 
\pagebreak
\begin{cor}
	Let $R\geq R_I\vee R_V$ and $u \in \CF_\text{loc}$ be a solution of equation of $(\ref{eqn:Poisson})$ in $W_R$. Let $\eta \in C_c^\infty(W_R)$ be a cutoff function. Then there exists a constant $C_1 > 0$ such that for all $\alpha \geq 1$, 
	\begin{align}\label{ine:itr}
		\norm{\eta u}_{W_R,\alpha \rho}^{2\alpha} \leq& \alpha ^2 C_1\norm{\lambda^{-1}}_{W_R,q}\norm{\Lambda}_{W_R,p}|W_R|^\frac{2}{d} \\
		&\times (\norm{\nabla \eta}_\infty^2 \norm{u}_{W_R, 2\alpha p^*}^{2\alpha} + \norm{\nabla f}_\infty^2 \norm{u}_{W_R,2\alpha p^*}).  \nonumber
	\end{align}
\end{cor}

The general idea of the proof of the following proposition is similar to that of \cite[Corollary 2.2]{CD} but we need to use $(\ref{ine:vol})$. 
\begin{prop}
	Take $R\geq R_I \vee R_V$ so that (\ref{ine:vol}) and $(\ref{ine:wiso})$ hold. We write $W(R)$ for $W_R$. Suppose that $u$ is a solution of (\ref{eqn:Poisson}) in $W(R)$, and assume that $|\nabla f|\leq c_f/R$. Then there exist $\kappa \in (1,\infty), \gamma \in (0,1]$ and $C_2=C_2(c_f)>0$ such that
	\begin{align}\label{ine:mx}
		\norm{u}_{W(\sigma'R), \infty}\leq C_2 \Biggl(\frac{1 \vee \norm{\lambda^{-1}}_{W(R),q} \norm{\Lambda}_{W(R),p} }{(\sigma - \sigma')^2} \Biggr)^\kappa \norm{u}_{W(\sigma R), \rho}^\gamma \vee \norm{u}_{W(\sigma R),\rho}
	\end{align}
	for any fixed $1/2\leq \sigma' < \sigma \leq 1$.
\end{prop}
\begin{proof}
	Throughout the proof, we use $C$ to denote a constant depending only on $d$, $p$, $q$ and $C_1$ and may change from line to line.
	We are going to apply the inequality (\ref{ine:itr}) iteratively. For fixed $1/2 \leq \sigma ' \leq \sigma \leq 1$, and $k \in \mathbb{N}$ define
	\[\sigma _k = \sigma ' +2^{-k+1}(\sigma - \sigma'). \]
	It is immediate that $\sigma _k-\sigma_{k+1} = 2^{-k}(\sigma - \sigma ')$ and that $\sigma _1 = \sigma$, furthermore $\sigma _k \downarrow \sigma '$. Let $p^*=p/(p-1)$. Then we have that $\rho > 2p^*$. Set $\alpha _k = (\rho/2p^*)^k$, $k\geq 1$. By definition, we have $\alpha _k\geq1$ for $k\geq 1$. Let $\tilde{\eta}_k\colon \mathbb{R}^d \to \mathbb{R}$ be a smooth function which is identically $1$ on $B(0,\sigma_{k+1}R)$ and vanishing on $\mathbb{R}^d\setminus B(0,\sigma_kR)$ and satisfies $\norm{\nabla \tilde{\eta}_k}\leq \frac{2^k}{(\sigma-\sigma')R}$. Define a function $\eta _k$ by the restriction of $\tilde{\eta}_k$ to $W$. 
	An application of (\ref{ine:itr}) and of the relation $\alpha_k\rho = 2\alpha _{k+1}p^*$ yields
	\begin{align*}
	&\norm{u}_{W(\sigma_{k+1}R),2\alpha_{k+1}p^*} \\
	&\leq \left (C \frac{2^{2k}\alpha _k^2|W(\sigma _kR)|^{2/d}}{(\sigma-\sigma')^2R^2}\norm{\lambda^{-1}}_{W(\sigma_k R),q}\norm{\Lambda}_{W(\sigma_k R),p} \right)^{1/(2\alpha_k)}\norm{u}_{W(\sigma_k R),2\alpha_kp^*}^{\gamma_k} \\
	&\leq \left (C \frac{2^{2k}\alpha _k^2}{(\sigma-\sigma')^2} \norm{\lambda^{-1}}_{W(R),q}\norm{\Lambda}_{W(R),p} \right)^{1/(2\alpha _k)}\norm{u}_{W(\sigma _kR),2\alpha_kp^*}^{\gamma_k},
	\end{align*}
	where $\gamma _k=1$ if $\norm{u}_{W(\sigma _kR),2\alpha_kp^*}\geq 1$ and $\gamma_k=1-1/\alpha_k$ otherwise. Iterating the above inequality and stop at $k=1$, we get
	\begin{align*}
	 \norm{u}_{W(\sigma _{j+1}R),2\alpha_{j+1}p^*}\leq \prod _{k=1}^j \left (C \frac{(\rho/p^*)^{2k}}{(\sigma-\sigma')^2}\norm{\lambda^{-1}}_{W(R),q}\norm{\Lambda}_{W(R),p} \right)^{1/(2\alpha _k)}\norm{u}_{W(\sigma R),\rho}^{\prod _{k=1}^j\gamma_k}. 
	\end{align*}
	Observe that $\kappa := \frac{1}{2}\sum \frac{1}{\alpha_k}< \infty, \; \sum \frac{k}{\alpha_k}<\infty$. Using (\ref{ine:vol}),  we have
	\begin{align*}
	&\norm{u}_{W(\sigma'R), 2\alpha_jp^*}\leq \left (\frac{|W(\sigma R)|}{|W(\sigma'R)|} \right )^{1/(2\alpha_jp^*)}\norm{u}_{W(\sigma R),2\alpha_jp^*}\\
	&\leq \left( \frac{C_V (\sigma R)^d}{(\sigma'R)^d} \right)^{1/(2\alpha_jp^*)} \norm{u}_{W(\sigma R),2\alpha_jp^*} \leq K\norm{u}_{W(\sigma R),2\alpha_jp^*}, 
	\end{align*}
	for some $K>0$ and all $j\geq 1$. Hence, taking the limit as $j \to \infty$ gives the inequality
	\begin{align*}
	\norm{u}_{W(\sigma'R),\infty} \leq C\left( \frac{1\vee \norm{\lambda^{-1}}_{W(R),q}\norm{\Lambda}_{W(R),p}}{(\sigma -\sigma')^2}  \right)^\kappa \norm{u}_{W(\sigma R),\rho}^{\prod_{k=1}^\infty\gamma_k},
	\end{align*}
	Define $\gamma = \prod _{k=1}^\infty(1-1/\alpha_k)\in (0,1]$, then $0 < \gamma \leq \prod _{k=1}^\infty \gamma_k$ and the above inequality can be written as
	\begin{align*}
	\norm{u}_{W(\sigma'R),\infty} \leq C\left( \frac{1\vee \norm{\lambda^{-1}}_{W(R),q}\norm{\Lambda}_{W(R),p}}{(\sigma -\sigma')^2}  \right)^\kappa \norm{u}_{W(\sigma R),\rho}^\gamma\vee\norm{u}_{W(\sigma R),\rho},
	\end{align*}
	which is the desired inequality.
\end{proof}

The main goal of this section is the following inequality. It is proved as in \cite[Corollary 2.2]{CD}.
\begin{cor}[maximal inequality]
	Let $R \geq R_I \vee R_V$. Suppose that $u$ is a solution of (\ref{eqn:Poisson}) in $W(R)$. Then, for all $\alpha \in (0,\infty)$ and for any $1/2 \leq \sigma'<\sigma<1$ there exist $C'=C'(p,q,d,c_f)>0,$ $\gamma'=\gamma'(\gamma,\alpha,\rho)$ and $\kappa' = \kappa'(\kappa,\alpha,\rho)$, such that 
	\begin{align}\label{ine:max}
	\norm{u}_{W(\sigma'R),\infty} \leq C'\left( \frac{1\vee \norm{\lambda^{-1}}_{W(R),q} \norm{\Lambda}_{W(R),p}}{(\sigma -\sigma')^2}  \right)^{\kappa'} \norm{u}_{W(\sigma R),\alpha}^{\gamma'}\vee\norm{u}_{W(\sigma R),\alpha}.
	\end{align}
\end{cor}

\section{Sublinearity}\label{sec:sublinear}
In this section, we will prove the sublinearity of the corrector and non-degeneracy of the covariance matrix. First, we will prove the $ L^p $-sublinearity. In \cite{CD}, the proof relies on  the fact that $ \E[U^k] = 0 $. However, as we mentioned in Remark \ref{rem:pot}, since we consider other measure $ \HP $, we have $ \HE[U^k] \neq 0$. To overcome this issue, we will take another approach.  

\begin{lem}\label{lem:supC}
	Let $\phi \in \mathcal{C}$ and define $ \Phi (x, \omega) $ by 
	\[ \Phi (x,\omega) = \int_0^1 |x \cdot \nabla \phi(\tau_{tx}\omega)|dt. \]
	Then,
	\[ \sup_{x \in \mathbb{R}^d}\Phi(x,\omega) < \infty \]
	holds $\HP$-almost surely $\omega$.
\end{lem}
\begin{proof}
	Let
	\[ \phi(\omega)=\int _{\mathbb{R}^d}f(\tau_x\omega)\eta(x)dx, \quad \eta \in C_c^\infty(\mathbb{R}^d) \]
	and set $ K = \text{supp}\;\eta $. Denote the diameter of $ K $ by $ \delta $. Then we have
	\begin{align*}
	\Phi(x,\omega)  &=  \int_0^1 |\langle x, \nabla \phi(\tau_{tx}\omega) \rangle |dt  \\
	&\leq  \int_0^1 \int_{\mathbb{R}^d}|\langle f(\tau_{tx+y}\omega)x, \nabla\eta(y)\rangle |dydt  \\
	&= \int_0^1 \int_{K}|\langle f(\tau_{tx+y}\omega)x, \nabla\eta(y) \rangle |dydt   \\
	&\leq \int_0^1\int_K \delta\norm{f}_\infty\norm{\nabla\eta}_\infty dydt \\
	&\leq \delta|K|\norm{f}_\infty\norm{\nabla\eta}_\infty.
	\end{align*}
	Hence we get the desired result.
\end{proof}

Set $\chi _\epsilon(x,\omega)=\epsilon \chi(x/\epsilon,\omega)$. To prove the $ L^p $-sublinearity, we recall a functional analysis result. Let $ \mathscr{B}_1, \mathscr{B}_2 $ be Banach spaces and $ T\colon \mathscr{B}_1 \to \mathscr{B}_2 $ a compact operator. Then, for each sequence $ \{x_n\}_n  \subset \mathscr{B}_1 $ such that $ x_n \to x $ weakly, we have that $ Tx_n \to Tx $  strongly.

\begin{lem}[$L^p$-sublinearity]\label{slt:p}
	For $\HP$-a.s. $\omega$ and $R \geq R_0(\omega)$, 
	\[\lim _{\epsilon \to 0}\norm{\chi _\epsilon^k(x,\omega)}_{2p^*,W_R}=0 .\]
\end{lem}
for $ k = 1,\dots, d $. 
\begin{proof}
	It is enough to show that for any $\eta \in C_c^\infty(W'_R)$ we have
	\begin{align}\label{eqn:weak}
	\lim_{\epsilon \to 0 }\frac{1}{|W'_R|}\int_{W'_R}\chi_\epsilon(x,\omega)\eta(x)dx =0.
	\end{align}
	
	Indeed, the above property implies the weak convergence $ \chi_\epsilon \to 0 $ in $L^2(W'_R)$. This gives the strong convergence in $ L^{2p^*}(W'_R) $, because $ W^{1,\frac{2q}{q+1}}(W'_R) $ is compactly embedded in $ L^{2p^*}(W'_R) $ and the sequence $ \{\chi_\epsilon\}_{\epsilon>0} $ is bounded in $ W^{1,\frac{2q}{q+1}}(W'_R) $. To show the equality (\ref{eqn:weak}), we will be apart from \cite{CD}. The following argument is motivated by \cite{DNS}.  
	
	We have
	\begin{align*}
	&\left| \frac{1}{|W'_R|} \int _{W'_R} \chi_\epsilon(x,\omega)\eta(x)dx  \right| \\
	&= \left| \frac{\epsilon^{d+1}}{|W'_R|} \int _{\frac{1}{\epsilon}W'_R} \chi(z,\omega)\eta(\epsilon z)dz  \right| \\
	&\leq \frac{\epsilon}{|\frac{1}{\epsilon}W'_R|} \int _{\frac{1}{\epsilon}W'_R} |\chi(z,\omega)\eta(\epsilon z)|dz \\
	&\leq \frac{\epsilon}{|\frac{1}{\epsilon}W'_R|} \int _{\frac{1}{\epsilon}W'_R}\int_0^1 | z\cdot	U^k(\tau_{tz}\omega)\eta(\epsilon z)|dtdz\\
	&\leq \norm{\eta}_\infty  \frac{\epsilon}{|\frac{1}{\epsilon}W'_R|}\int_0^1\int _{\frac{1}{\epsilon}W'_R} |z \cdot U^k(\tau_{tz}\omega)|dzdt.
	\end{align*} 
	Now, let $ \phi _n \in \mathcal{C} $ be an approximate sequence such that $ \nabla \phi _n  \to U^k $  in $ \hl $ and set 
	\[ \Phi_n(x,\omega) = \int _0^1 |x\cdot \phi_n(\tau_{tx}\omega)|dt. \] 
	Then we compute 
	\begin{align*}
	&\frac{\epsilon}{|\frac{1}{\epsilon}W'_R|}\int_0^1\int _{\frac{1}{\epsilon}W'_R} |z \cdot U^k(\tau_{tz}\omega)|dzdt\\
	&\leq \frac{\epsilon}{|\frac{1}{\epsilon}W'_R|}\int_0^1\int _{\frac{1}{\epsilon}W'_R} |z \cdot \nabla\phi_n(\tau_{tz}\omega)|dzdt \\
	&+\frac{\epsilon}{|\frac{1}{\epsilon}W'_R|}\int_0^1\int _{\frac{1}{\epsilon}W'_R} |z\cdot U^k(\tau_{tz}\omega)-z \cdot \nabla\phi_n(\tau_{tz}\omega)|dzdt.  
	\end{align*}
	The first term of the right hand side is bounded above by
	\begin{align*}
	&\frac{\epsilon}{|\frac{1}{\epsilon}W'_R|}\int_0^1\int _{\frac{1}{\epsilon}W'_R} |z \cdot \nabla\phi_n(\tau_{tz}\omega)|dzdt\\
	&= \frac{\epsilon}{|\frac{1}{\epsilon}W'_R|}\int _{\frac{1}{\epsilon}W'_R}\int_0^1 |z \cdot \nabla\phi_n(\tau_{tz}\omega)|dtdz \\
	&= \frac{\epsilon}{|\frac{1}{\epsilon}W'_R|}\int _{\frac{1}{\epsilon}W'_R}\Phi_n(z,\omega)dtdz \\
	&\leq \epsilon \sup_{z\in \mathbb{R}^d}\Phi_n(z,\omega).
	\end{align*}
	Since $ \sup_{z\in \mathbb{R}^d}|\Phi_n(z,\omega)| < \infty $ by Lemma \ref{lem:supC}, this term tends to zero as $ \epsilon \to 0 $. \\
	For the second term, using the Cauchy-Schwarz inequality, 
	we estimate
	\begin{align*}
	&\frac{\epsilon}{|\frac{1}{\epsilon}W'_R|}\int_0^1\int _{\frac{1}{\epsilon}W'_R} |z\cdot U^k(\tau_{tz}\omega)-z \cdot \nabla\phi_n(\tau_{tz}\omega)|dzdt. \\
	&=\frac{1}{|\frac{1}{\epsilon}W'_R|}\int_0^1\int _{\frac{1}{\epsilon}W'_R} |(\epsilon z)\cdot (U^k(\tau_{tz}\omega)- \nabla\phi_n(\tau_{tz}\omega))|dzdt. \\
	&\leq \frac{1}{|\frac{1}{\epsilon}W'_R|}\int_0^1\int _{\frac{1}{\epsilon}W'_R} \sqrt{|\epsilon z||U^k(\tau_{tz}\omega)- \nabla\phi_n(\tau_{tz}\omega)|}dzdt.
	\end{align*}
	Using the fact that $|z| \leq R/\epsilon$ for $z \in \frac{1}{\epsilon} W_R$ and the change of variables, we have
	\begin{align*}
	&\frac{1}{|\frac{1}{\epsilon}W'_R|}\int_0^1\int _{\frac{1}{\epsilon}W'_R} \sqrt{|\epsilon z||U^k(\tau_{tz}\omega)- \nabla\phi_n(\tau_{tz}\omega)|}dzdt \\
	&\leq R^\frac{1}{2}\frac{1}{|\frac{1}{\epsilon}W'_R|}\int_0^1\int _{\frac{1}{\epsilon}W'_R} \sqrt{|U^k(\tau_{tz}\omega)- \nabla\phi_n(\tau_{tz}\omega)|}dzdt \\
	&= R^\frac{1}{2}\int_0^1 \frac{1}{|\frac{t}{\epsilon}W'_R|}\int _{\frac{t}{\epsilon}W'_R} \sqrt{|U^k(\tau_{u}\omega)-\nabla\phi_n(\tau_{u}\omega)|}dudt.
	\end{align*}
	By the volume regularity ((1) of Assumption \ref{asm:volIso}), this is bounded above by 
	\begin{align*}
	&C R^\frac{1}{2}\int_0^1 \frac{1}{|B(0,\frac{tR}{\epsilon})|}\int _{B(0,\frac{tR}{\epsilon})} \sqrt{|U^k(\tau_{u}\omega)-\nabla\phi_n(\tau_{u}\omega)|}\1{0 \in	W'(\tau_u\omega)}dudt.
	\end{align*} 
	By the ergodic theorem, we have
	\[ \lim_{\epsilon \to 0 }\frac{1}{|B(0,\frac{tR}{\epsilon})|}\int _{B(0,\frac{tR}{\epsilon})} \sqrt{|U^k(\tau_{u}\omega)-\nabla\phi_n(\tau_{u}\omega)|}\1{0 \in	W'(\tau_u\omega)}du = \HE[\sqrt{|U^k-\nabla \phi_n|}] .\]
	Letting $n$ tend to infinity, this term converges to zero. Now we have
	\begin{align*}
	&\left| \frac{1}{|W'_R|} \int _{\mathbb{R}^d} \chi_\epsilon(x,\omega)\eta(x)dx  \right| \\
	&\leq \norm{\eta}_\infty \Big( \epsilon \sup_{z\in \mathbb{R}^d}|\Phi_n(z,\omega)| \\
	&+ C_V^{-1} R\int_0^1 \frac{1}{|B(0,\frac{tR}{\epsilon})|}\int _{B(0,\frac{tR}{\epsilon})} |U^k(\tau_{u}\omega)-\nabla\phi_n(\tau_{u}\omega)|\1{0 \in	W'(\tau_u\omega)}dudt \Big).
	\end{align*}
	First let $\epsilon$ tend to zero and then $ n $ tend to infinity, we get (\ref{eqn:weak}) and obtain the result.
\end{proof}

\begin{prop}[$L^\infty$-sublinearity]\label{prop:slti}
	For all $R>0$,
	\[\lim _{\epsilon \to 0}\sup _{x \in W'_R}|\chi _\epsilon (x,\omega)|=0,\quad \HP \text{-a.s.}\]
\end{prop}
\begin{proof}
	Observe that by Proposition \ref{prop:har},   the function $\chi_\epsilon^k(x,\omega)$ is a solution of
	\[\CE^\omega (\chi_\epsilon, \phi)=\int _{W'} \langle a(\tau_x\omega)\nabla f_k, \nabla \phi \rangle\; dx \]
	in $W'_R$ for all $\epsilon >0$. Here $f_k(x)=x_k$ and $\displaystyle\phi \in H^1(W'_R)\cap C_c^\infty (\overline{W_R})$. We first consider the case $R \geq R_0(\omega)$. Since $|\nabla f_k|\leq 1$, by Lemma \ref{slt:p} we have
	\[\lim _{\epsilon \to 0} \norm{\chi_\epsilon^k(x,\omega)}_{2p^*,W'_R}=0.\]
	Therefore, for $R \geq R_0(\omega)$, we get the desired result by the maximal inequality (\ref{ine:max}) with $\alpha = 2p^*$. It remains to treat the case that $R \in (0,R_0(\omega))$, but in that case, the desired bound immediately follows from the fact that $\sup _{x \in W'_R}|\chi _\epsilon (x,\omega)| \leq \sup _{x \in W'_{R_0}}|\chi _\epsilon (x,\omega)|$ for $R \in (0, R_0)$.
\end{proof}

To prove the positive-definiteness of the covariance matrix, we show the following lemma. Recall that $T_x$ is the operator on $L^\infty(\HP)$ defined by $T_xG(\omega) = G(\tau_x\omega)$.

\begin{lem}\label{lem:vanishPathInt}
	Let $G\colon \Omega \to \mathbb{R}$ be an integrable random variable such that $G(\omega) = 0$ for $\HP$-almost all $\omega \in \OH$. 
	Fix $x \in \mathbb{R}^d$. 
	Then for $\HP$-almost all $\omega$, 
	\begin{align*}
		\int_0^1 (T_{\gamma(t)}G) \1{\gamma(t) \in W'(\omega)}\gamma'(t) dt = 0,
	\end{align*}
	holds for all smooth path $\gamma\colon [0,1] \to \mathbb{R}^d$ satisfying $\gamma(0)=0$ and $\gamma(1)=x$. 
\end{lem}
\begin{proof}
	First we fix a smooth path $\gamma$. Then for all $F \in L^\infty(\HP)$, by Fubini's theorem and the definition of $\HP$, we have
	\begin{align*}
		&\HE \biggl[F \int_0^1  (T_{\gamma(t)}G)\1{\gamma(t)\in W'} \gamma'(t)  dt \biggr] \\
		&= \int_0^1\HE [F  (T_{\gamma(t)}G)\1{\gamma(t)\in W'} \gamma'(t) ]dt  \\
		&= \int_0^1\E [F  (T_{\gamma(t)}G)\1{\gamma(t)\in W'} \gamma'(t) \1{0\in W'}]\P(0\in W')^{-1}dt.
	\end{align*}

	Using the translation invariance of $\P$, the last line equals to 
	\begin{align*}
		&\int_0^1\E [(T_{-\gamma(t)}F)G\1{0\in W'}\1{-\gamma(t)\in W'} \gamma'(t) ]\P(0\in W')^{-1}dt\\
		&= \int_0^1\HE \bigl[(T_{-\gamma(t)}F)  G \gamma'(t) \1{-\gamma(t)\in W'}\bigr]dt\\
		&= \int_0^1\HE[(T_{-\gamma(t)}F)\cdot 0 \cdot \gamma'(t)\1{-\gamma(t)\in W'}]dt\\
		&=0.
	\end{align*}

	Since $F \in L^\infty(\HP)$ is arbitrary, we obtain 
	\begin{align*}
		\int_0^1 \langle T_{\gamma(t)}G, \gamma'(t)\rangle \1{\gamma(t) \in W'}dt = 0,
	\end{align*}
	for $\HP$-a.s. $\omega$.
	This implies the result because the collection of all smooth path is separable with respect to the supremum norm (that is, for a path $\gamma\colon [0,1] \to \mathbb{R}^d$, define its norm by $\sup_{t\in[0,1]}|\gamma (t)|$).
\end{proof}
	
\begin{prop}
	The covariance matrix $ \mathbb{D} $ is positive-definite.
\end{prop}

\begin{proof}
The general idea of the proof goes back to \cite{DNS}.
Assume that there is a vector $v \in \mathbb{R}^d$ with $|v|=1$ such that $\langle v,\mathbb{D}v \rangle=0$. We have
\begin{align*}
\langle v , \mathbb{D}v \rangle &= \sum _{i=1}^d v_i(\mathbb{D}v)_i \\
&= \sum _{i=1}^d v_i\left(\sum _{j=1}^d d_{ij}v_j \right) \\
&= \sum _{k=1} \HE [\langle a\xi, \xi \rangle],
\end{align*} 
where $\xi = (\langle v, D_1y\rangle, \dots ,\langle v, D_d y\rangle)$,  $y = (y^1, \dots , y^d)$ and $D_ky = (D_ky^1, \dots, D_ky^d)$. 
Thus, we have $\langle v, D_ky(0,\omega)\rangle=0$ for $k=1,\dots,d$ $\HP$-a.s. 
Take $\omega \in \OH$ and smooth path $\gamma_\omega \colon [0,1] \to \mathbb{R}^d$ such that $\gamma_\omega(0)=0$, $\gamma_\omega(1)=x$ and $\gamma_\omega([0,1])\subset W'(\omega)$. We denote the $k$-th component of $\gamma_\omega'$ by $\gamma_{\omega,k}'$. Since $D_ky(x,\omega)=D_ky(0,\tau _x\omega)$ and
\begin{align*}
	 y^i(x,\omega) = \int_0^1 \langle \nabla y^i(\gamma_\omega(t),\omega),\gamma_\omega'(t)\rangle dt = \sum_{k=1}^d\int_0^1D_k y^i(0,\tau_{\gamma_\omega(t)}\omega)\gamma_{\omega,k}'(t)\1{\gamma_\omega(t)\in W'(\omega)}dt,
\end{align*}
it holds that 
\begin{align*}
	\langle v, y(x,\omega)\rangle &= \sum_{i=1}^d\sum_{k=1}^d\int_0^1 v_iD_k y^i(0,\tau_{\gamma(t)}\omega)\gamma_{\omega,k}'(t)\1{\gamma_\omega(t)\in W'(\omega)}dt\\
	&= \sum_{k=1}^d \int_0^1 \langle v, D_ky(0,\tau_{\gamma_\omega(t)}\omega)\rangle \1{\gamma_\omega(t)\in W'}\gamma_\omega'(t)dt.
\end{align*}
Thus it follows from Lemma \ref{lem:vanishPathInt} that we have $\langle v, y(x,\omega)\rangle = 0$.
Hence $\langle v, x\rangle = \langle v, \chi(x,\omega)\rangle$ $\HP$-a.s. However, it implies 
\[ 0 = \lim _{\epsilon \to 0 }\sup _{x \in W'_R }\langle v, \epsilon \chi  (x/\epsilon, \omega)\rangle  = \lim _{\epsilon \to 0 }\sup _{x \in W'_R }\langle v, \epsilon \left(x/\epsilon\right) \rangle =\lim _{\epsilon \to 0 }\sup _{x \in W'_R }\langle v, x \rangle  >0 ,\] a contradiction.
\end{proof}

\section{Proof of the main theorem}
	Recall that we have the decomposition  
\[ \epsilon X_{\cdot/\epsilon^2}^\omega = \epsilon\chi(X_{\cdot/\epsilon^2}^\omega,\omega)+\epsilon y(X_{\cdot/\epsilon^2}^\omega,\omega). \]
First, we show that the martingale $ \epsilon y(X_{\cdot/\epsilon^2}^\omega,\omega) $ converges to a Brownian motion with covariance matrix $\mathbb{D}$. 

\begin{lem}\label{lem:martconv}
	Set $M_t^\omega = y(X_t^\omega, \omega) $. Then $\epsilon M_{\cdot/\epsilon^2}^\omega$ converges to a Brownian motion with the covariance matrix $\mathbb{D}$.
\end{lem}
\begin{proof}
	By Corollary \ref{cor:mart} and the ergodic theorem for the environment process $ \{\tau_{X_t^\omega}\omega\}_t $ (see \cite[Proposition 2.1]{O}), we have
	\begin{align*}
		\langle M_\cdot^{k,\epsilon}, M_\cdot^{\ell,\epsilon} \rangle_t &= \epsilon^2\int _0^{t/\epsilon^2} \sum _{i=1}^d\langle a(\tau_{X_s^\omega}\omega)\nabla y^k(X_s,\omega),\nabla y^\ell(X_s,\omega)\rangle ds\\
		&\to \HE [\langle a(\omega)\nabla y^k(0,\omega), \nabla y^\ell(0,\omega)\rangle].
	\end{align*}
	Hence, we get the result by \cite[Theorem 5.1.]{H}.
\end{proof}

It remains to show that the corrector $\epsilon \chi (X_{t/\epsilon^2}^\omega , \omega)$ converges to zero in distribution. For that the sublinearity of the corrector will play a major role.

Let $T>0$ be a fixed time horizon. We claim that for all $\delta >0$
\begin{align}\label{eqn:cor}
	\lim _{\epsilon\to 0} P_0^\omega\left(\sup_{0 \leq t \leq T} |\epsilon\chi(X_{t/\epsilon^2}^\omega ,\omega) | > \delta \right )= 0.
\end{align}
Denote by $\tau _R^{\epsilon,\omega}$ the exit time of $\epsilon X_{t/\epsilon^2}^\omega$ from  $W'_R$. Observe that 
\begin{multline}
	\limsup _{\epsilon\to 0} P_0^\omega \left (\sup_{0 \leq t \leq T}|\epsilon \chi (X_{t/\epsilon^2}^\omega,\omega)| > \delta \right)  \\
	\leq \limsup_{\epsilon \to 0}P_0^\omega \left (\bigg|\sup_{0 \leq t \leq \tau _R^{\epsilon,\omega}}\epsilon \chi (X_{t/\epsilon^2}^\omega,\omega)\bigg| > \delta \right) + \limsup _{\epsilon \to 0} P_0^\omega \left(\bigg|\sup_{0\leq t \leq T}\epsilon X_{t/\epsilon^2}^\omega\bigg| > R \right) \label{ine:dif}.
\end{multline}

By Proposition \ref{prop:slti}, we have
\[\lim_{\epsilon \to 0 }\sup_{0 \leq t \leq \tau _R^{\epsilon,\omega}}|\epsilon\chi(X_{t/\epsilon^2}^\omega,\omega)| = 0 \]
and therefore $ \HP $-almost surely
\[ \limsup_{\epsilon\to 0} P_0^\omega\left(\sup_{0 \leq t \leq \tau _R^{\epsilon,\omega}} |\epsilon\chi(X_{t/\epsilon^2}^\omega,\omega)|>\delta \right)=0.\]

To obtain a bound of the second term of $(\ref{ine:dif})$,  we use again Proposition \ref{prop:slti} to say that there  exists $ \tilde{\epsilon}(\omega)>0 $, which may depend on $ \omega $ such that for all $ \epsilon < \tilde{\epsilon} $ we have $ \sup_{0 \leq t \leq \tau _R^{\epsilon,\omega}} |\epsilon\chi(X_{t/\epsilon^2}^\omega,\omega) |<1 $. For such $ \epsilon $ we have $ \HP $-almost surely 
\begin{align*}
	P_0^\omega\left(\sup_{0 \leq t \leq T} |\epsilon X_{t/\epsilon^2}^\omega|\geq R\right) &= P_0^\omega(\tau_R^{\epsilon,\omega} \leq T) \\
	&= P_0^\omega \left(\tau_R^{\epsilon,\omega} \leq T,\sup_{0 \leq t \leq \tau _R^{\epsilon,\omega}}|\epsilon y(X_{t/\epsilon^2}^\omega,\omega)|> R-1\right) \\
	&\leq P_0^\omega \left(\sup_{0 \leq t \leq T}|\epsilon y(X_{t/\epsilon^2}^\omega,\omega)|> R-1 \right).
\end{align*}

Thanks to Lemma \ref{lem:martconv}, the process $ \epsilon y(X_{\cdot/\epsilon^2}^\omega,\omega) $ converges in distribution under $ P_0^\omega $ to a non-degenerate Brownian motion with the deterministic covariance matrix given by $ \mathbb{D} $. Hence there exist positive constants $ c_1,c_2 $ independent of $ \epsilon $ and $ \omega $ such that  
\[ \limsup_{\epsilon \to 0} P_0^\omega \left (\sup_{0\leq t \leq T}|\epsilon y(X_{t/\epsilon^2}^\omega,\omega)|>R-1 \right ) \leq c_1e^{-c_2R} ,\]
from which it follows
\[ \limsup_{\epsilon \to 0}P_0^\omega\left (\sup_{0 \leq t \leq T}|\epsilon X_{t/\epsilon^2}^\omega|>R\right ) \leq c_1e^{-c_2R} .\]
Therefore, 
\[ \limsup_{\epsilon \to 0}P_0^\omega\left (\sup_{0 \leq t \leq T}|\epsilon \chi (X_{t/\epsilon^2}^\omega,\omega)|>\delta \right ) \leq c_1e^{-c_2R} \]
and since $ R>1 $ was arbitrary, the claim (\ref{eqn:cor}) follows, namely the corrector converges to zero in law under $ P_0^\omega $, $ \HP $-almost surely. 

The convergence to zero in law of the corrector $ \epsilon\chi(X_{\cdot/\epsilon^2}^\omega,\omega)$, combined with the fact that $ \epsilon y(X_{\cdot/\epsilon^2}^\omega,\omega) $ satisfies an invariance principle $ \HP $-almost surely and that $ \epsilon X_{\cdot/\epsilon^2}^\omega = \epsilon\chi(X_{\cdot/\epsilon^2}^\omega,\omega)+\epsilon y(X_{\cdot/\epsilon^2}^\omega,\omega) $, implies that the family $ \epsilon X_{\cdot/\epsilon^2}^\omega $ under $ P_0^\omega $ satisfies an invariance principle $ \HP $-almost surely with the same limiting law.

$ $\\
\textbf{Acknowledgements} The author would like to thank Professor Hideki Tanemura for suggesting the problem and fruitful discussions. This work was supported by Keio University Doctoral Student Get-in-Aid Program. 
\begin{flushleft}
\end{flushleft}
	
\bibliographystyle{plain}
\bibliography{QIPreference}

\end{document}